\documentclass[reqno]{amsart}

\usepackage[dvipsnames]{xcolor}
\usepackage{amsmath}
\usepackage{enumerate}
\usepackage{enumitem}
\usepackage{amssymb}
\usepackage{extarrows}
\usepackage{amsthm}
\usepackage{thmtools}
\usepackage{mathrsfs}
\usepackage[colorlinks, linkcolor = SeaGreen, citecolor = SeaGreen]{hyperref}
\usepackage{bm}
\usepackage{tikz-cd}
\hypersetup{colorlinks = true}

\newtheorem{theorem}{Theorem}[section]

\newtheorem{lemma}[theorem]{Lemma}

\newtheorem{proposition}[theorem]{Proposition}
\newtheorem{question}[theorem]{Question}
\theoremstyle{definition}

\theoremstyle{remark}
\newtheorem{remark}[theorem]{Remark}

\newcommand{\B}{\mathbb B_n}

\newcommand{\Hh}{\mathcal H}

\newcommand{\ip}[2]{\left\langle #1,#2\right\rangle}

\numberwithin{equation}{section}

\title[Zero products of Toeplitz operators]{Zero products of Toeplitz operators with harmonic symbols on the Bergman space of the unit ball}

\author[Wang]{Jiawei Wang}
\address{School of Mathematical Sciences, Inner Mongolia University, Hohhot, 010021, P. R. China}
\email{jiaweiwang1516@163.com}

\author[Zhao]{Xianfeng Zhao}
\address{College of Mathematics and Statistics, Chongqing University, Chongqing, 401331, P. R. China}
\email{xianfengzhao@cqu.edu.cn}


\makeatletter

\newcommand{\Rmnum}[1]{\expandafter\@slowromancap\romannumeral #1@}
\@namedef{subjclassname@2020}{\textup{2020} Mathematics Subject Classification}
\makeatother
\date{\today}

\begin{document}
\begin{abstract}
In this paper, using local boundary factorizations of harmonic functions and singular-value estimates for Toeplitz products, we study the zero-product problem for Toeplitz operators on the Bergman space of the unit ball in $\mathbb{C}^n$. We prove that a finite product of such operators is zero only if one of the symbols is identically zero, provided that all the symbols are harmonic on the open unit ball and continuous on the closed unit ball.
\end{abstract}
\maketitle

\tableofcontents

\section{Introduction}

Let $n$ be a positive integer and $\B=\{z\in\mathbb C^n:|z|<1\}$ be the open unit ball in $\mathbb C^n$ with normalized volume measure $dV$.
For $z=(z_1, z_2, \cdots,z_n)\in \mathbb C^n$ and $w=(w_1, w_2,\cdots,w_n)\in \mathbb C^n$, the inner product on
$\mathbb C^n$ is given by
\[
\ip zw=\sum_{j=1}^n z_j\overline{w}_j.
\]
Let $\mathbb S_n$ be the unit sphere in $\mathbb C^n$ and $d\sigma$ be the normalized surface measure. Noting that $d\sigma$ is the rotation-invariant surface measure on $\mathbb S_n$, we have the following polar integration formula:
\begin{equation}\label{polar}
\int_{\B}F(z)\,dV(z)
=2n\int_0^1 t^{2n-1}\int_{\mathbb S_n}F(t\zeta)\,d\sigma(\zeta)\,dt,
\end{equation}
where $F$ is integrable on $\mathbb B_n$ with respect to  $dV$, i.e., $F\in L^1(\mathbb B_n, dV)$.

The Hilbert space $L^2(\mathbb{B}_n)=L^2(\mathbb{B}_n,dV)$ consists of all
square-integrable functions on $\mathbb{B}_n$ with respect to $dV$,
equipped with the inner product and norm
\[
\langle f,g\rangle
=\int_{\mathbb{B}_n} f(z)\overline{g(z)}\,dV(z),
\ \ \ \
\|f\|_2
=\left(\int_{\mathbb{B}_n}|f(z)|^2\,dV(z)\right)^{\frac{1}{2}}.
\]
 The Bergman space of the unit ball is defined by
\[
L_a^2(\mathbb B_n)=\{f\in L^2(\B): f\text{ is holomorphic on }\B\},
\]
which is a reproducing kernel  Hilbert space, and the reproducing kernel at $w\in\B$ is given by
\begin{equation*}\label{eq:kernel}
K_w(z)=\frac{1}{(1-\ip zw)^{n+1}},\  \ \ \ z\in \mathbb B_n.
\end{equation*}
Moreover, the normalized reproducing kernel  is
\begin{equation*}\label{eq:normalized-kernel}
k_w(z)=\frac{K_w(z)}{\|K_w\|_2}
=\frac{(1-|w|^2)^{\frac{n+1}{2}}}{(1-\langle z,w\rangle)^{n+1}},\  \ \ \ z\in \mathbb B_n.
\end{equation*}

Let $P$ be the orthogonal projection  from $L^2(\B)$ onto $L_a^2(\mathbb B_n)$. Using the reproducing kernel introduced above, $P$ can be represented  by
$$
(Pf)(z)=\int_{\B}\frac{f(w)}{(1-\ip zw)^{n+1}}\,dV(w),\ \ \ f\in L^2(\B).
$$
Let $L^\infty(\B)= L^\infty(\B, dV)$ denote the Banach algebra of essentially bounded functions on $\B$ with respect to the measure $dV$. For $\varphi\in L^\infty(\B)$, the \emph{Toeplitz operator} and the \emph{Hankel operator} with symbol $\varphi$ on $L_a^2(\mathbb B_n)$ are defined by
$$T_{\varphi}f=P(\varphi f)$$
and
$$H_{\varphi}f=(I-P)(\varphi f).$$
Clearly, $\|T_\varphi\|\leqslant \|\varphi\|_\infty$ and $\|H_\varphi\|\leqslant \|\varphi\|_\infty$.

For $\varphi$ and $\psi$ in $L^\infty(\mathbb B_n)$, the operators $T_{\varphi}T_{\psi}-T_{\psi}T_{\varphi}$ and $T_{\varphi}T_{\psi}-T_{\varphi \psi}$ are called the commutator and the semi-commutator of
$T_{\varphi}$ and $T_{\psi}$, respectively.
Analogously to the case of the Bergman space of the open unit disk $\mathbb D$,
by using the decomposition of  the multiplication operators on
\[
L^2(\mathbb{B}_n)
=
L_a^2(\mathbb{B}_n)
\oplus
\bigl(L_a^2(\mathbb{B}_n)\bigr)^\perp,
\]
we obtain the following:

\begin{equation}\label{eq:Hankel identity}
  T_{\varphi \psi}-T_{\varphi}T_{\psi}=H_{\overline \varphi}^{\,*}H_{\psi}.
\end{equation}
Moreover, it follows from  \cite[Theorem 1]{Coburn1973}  that the commutator $T_{\varphi} T_{\psi}-T_{\psi}T_{\varphi}$ and the semi-commutator $T_{\varphi}T_{\psi}-T_{\varphi \psi}$ are both compact on $L_a^2(\mathbb B_n)$ if $\varphi$ and $\psi$ are continuous on the closed
unit ball $\overline{\mathbb B}_n$ (i.e., $\varphi, \psi\in C(\overline{\mathbb B}_n)$). Therefore, the Hankel operator $H_{\varphi}$ is compact   when $\varphi\in C(\overline{\mathbb B}_n)$.

For further details about the Bergman space of the disk or the unit ball, as well as Toeplitz and Hankel operators acting on these spaces, we refer the reader to \cite{Coburn1973, Zhu, Zhu-book}.

In this paper, we investigate the \emph{zero-product problem} for Toeplitz operators
on the Bergman space \(L_a^2(\mathbb{B}_n)\). More precisely, we consider the following question:
\begin{question}\label{q:zero-product}
If \(\varphi_1, \varphi_2, \cdots,\varphi_m\in L^\infty(\B)\) satisfy
\[
T_{\varphi_1} T_{\varphi_2}\cdots T_{\varphi_m}=0,
\]
then must one of the functions \(\varphi_j\) be zero almost everywhere on \(\B\)?
\end{question}

The zero-product problem for Toeplitz operators originated in \cite{BH}.
In that paper, Brown and Halmos  established that if the product of two Toeplitz operators on the Hardy space of the unit disk is zero, then one of the symbols is zero almost everywhere. This leads to the following natural question, which has been of great interest to researchers:

\begin{question}\label{q:zero-product-Hardy}
If a product of $m$ Toeplitz operators on the Hardy space is equal to zero, must one of these operators be zero?
\end{question}

Question \ref{q:zero-product-Hardy} was answered in the affirmative for three operators by Axler in the 1970s, but his method becomes considerably more complicated when dealing with more operators.  More than 20 years later, using some function-theoretic techniques on the Hardy space, Guo \cite{Guo} studied the kernel spaces of the products of Toeplitz operators, thereby resolving  Question \ref{q:zero-product-Hardy} for the case  of $m=5$.  In 2000, Gu \cite{Gu} succeeded in resolving the  case of $m=6$  by further developing the techniques originally due to Guo \cite{Guo}. In 2009, this question was completely settled by Aleman and Vukoti\'{c} \cite{AV}, who employed techniques involving nearly-invariant subspaces and vector-valued Hardy spaces.

The zero-product problem for Toeplitz operators on the Bergman space
is also quite difficult. Even in the case \(\mathbb B_1=\mathbb D\),
the open unit disk, the answer to Question \ref{q:zero-product} remains unknown.
Since the harmonic extension gives a natural correspondence between the bounded functions on the unit circle and the bounded harmonic functions on the unit disk \cite{Dur}, Toeplitz operators on the Bergman space with bounded harmonic symbols behave quite similarly to those on the Hardy space. This makes harmonic symbols a natural class for studying the Bergman-space analogues of the zero-product problem for Toeplitz operators on the Hardy space.

Let us first recall the case of the Bergman space on the disk. In 2001, Ahern and \v{C}u\v{c}kovi\'{c} \cite{AC} solved the zero-product problem for two Toeplitz operators on the Bergman space with  bounded harmonic symbols via the Berezin transform. Subsequently, \v{C}u\v{c}kovi\'{c} \cite{Cu} showed that if $\Psi(z)=z^j-\overline{z}^\ell$ and $\Phi$ is an arbitrary bounded function such that $T_{\Phi}T_{\Psi}=0$ on the Bergman space, then $\Phi=0$. This result was later generalized to the case when $\Psi$ is an arbitrary nonzero polynomial in $z$ and $\overline{z}$, see \cite{Le2010} and \cite{LRY}.  Moreover, Ahern and \v{C}u\v{c}kovi\'{c} \cite{AZ1} also resolved the zero-product problem for two Toeplitz operators on the Bergman space when one of the symbols is assumed to be radial.

On the Bergman space of the unit ball in \(\mathbb{C}^n\), Choe and Koo \cite{CK} solved the zero-product problem for two Toeplitz operators with harmonic symbols that have continuous extensions to some part of the boundary of $\mathbb B_n$. Furthermore, in the case where symbols have Lipschitz continuous extensions to the boundary, they solved  the zero-product problem for multiple products with the number of factors depending on the dimension \(n\) of the underlying space. In addition, further partial results on the zero-product problem of Toeplitz operators can be found in \cite{CKL1, CKL2, CLNZ, Le2009, Das, CHS} for the Bergman space of the polydisk, the unit ball, an annulus and Reinhardt domains, and in \cite{Bo, Qin} for the Fock space.

In this paper, using local boundary factorizations of harmonic functions and singular-value estimates for Toeplitz products, we solve the zero-product problem for Toeplitz operators with harmonic symbols on the Bergman space $L_a^2(\mathbb B_n)$, under the assumption that all  symbols are continuous on the closed  unit ball. The main result of our paper is as follows.

\begin{theorem}\label{thm:main}
Suppose
$u_1, u_2, \cdots,u_m$ are complex-valued functions
that are harmonic on $\mathbb B_n$ and continuous on the closed unit ball
$\overline{\mathbb B}_n$.
If $$T_{u_1}T_{u_2}\cdots T_{u_m}=0$$
on the Bergman space $L_a^2(\mathbb B_n)$, then $u_j\equiv0$ on $\mathbb B_n$
for some $j\in\{1,2,\cdots,m\}$.
\end{theorem}

It is worth noting that, compared with \cite{CK}, the above theorem  removes the restriction on the number of factors and weakens the global Lipschitz assumption to continuity. Moreover, our result is new even in the case of the  Bergman space of the unit disk.

In the remainder of this section, we briefly describe the motivation for studying the zero-product problem for Toeplitz operators with harmonic symbols, and outline the approach and key ideas of the proof.

Assume  that  $u_1, u_2, \cdots,u_m$ are harmonic on $\mathbb B_n$ and continuous on
$\overline{\mathbb B}_n$. Moreover, suppose
$$
T_{u_1}T_{u_2}\cdots T_{u_m}=0
$$
on the Bergman space $L_a^2(\mathbb B_n)$.
Under this assumption, by using Coburn's theorem (\cite[Theorem 1]{Coburn1973}) we first obtain that
$$
u_1(\zeta)u_2(\zeta)\cdots u_m(\zeta)=0,
\qquad \zeta\in\mathbb S_n.
$$

Therefore, the next natural question is  what form a
harmonic function takes near an open portion of the
boundary on which it vanishes. Indeed, the answer comes from the Poisson formula:
$$
u(z)
=
(1-|z|^2)
\int_{\mathbb S_n}
\frac{u(\xi)}{|z-\xi|^{2n}}\,d\sigma(\xi), \ \ \ \ z\in \mathbb B_n.
$$
If $u=0$ on an open subset $W$ of the sphere $\mathbb S_n$, then we get that
$$
u(z)
=
(1-|z|^2)
\int_{\mathbb S_n\backslash W}
\frac{u(\xi)}{|z-\xi|^{2n}}\,d\sigma(\xi).
$$
Choose a smaller open subset of the sphere whose
closure is contained in $W$.
In a sufficiently small neighborhood  of this closure
in $\mathbb C^n$ (denoted by $U$), the variable $z$ remains a positive
distance from the region of integration. Thus,
the integral  defines a smooth function
$\phi$.
In the intersection of $U$ with the
ball $\mathbb B_n$, we have
$$
u(z)=\rho(z)\phi(z),
$$
where $\rho(z)=1-|z|^2$.
Observe that the function $\rho$ describes the
vanishing of these harmonic symbols near the boundary.

Returning to Toeplitz operators, the original symbol
locally has the form $\rho\phi$.
To use this factorization on $\mathbb B_n$,
we first extend the function $\phi$ obtained above and then
localize the operator product using a cutoff.
The factor we consequently need to study is $
T_{\rho\phi}$.
The localization result ensures that modifying the
symbols outside the chosen neighborhood, after
multiplication on the right by the local cutoff
Toeplitz operator, produces an operator difference
with rapidly decaying singular values.
To analyze how this factor affects the whole
product, we first consider the simplest case $\phi\equiv 1$, i.e., $T_{\rho\phi}=T_\rho$.
Then we  prove that, for a continuous function
$\phi$, we can write
$$
T_{\rho\phi}
=
T_\phi T_\rho+E,
\qquad
s_k(E)=o(k^{-\frac{1}{n}}) \ \ \ \  (k\rightarrow \infty).
$$
Here, $s_k(E)$ denotes the $k$-th singular value of
the compact operator $E$, arranged in nonincreasing
order.
Thus, we transform a factor $\rho$ in the symbol into a factor
$T_\rho$ in the operator product. Note that the properties of $T_\rho$  are relatively simple, as its symbol is radial. Therefore, the natural line of reasoning is as follows:
$$
\begin{array}{ccc}
T_{u_1}T_{u_2}\cdots T_{u_m}=0
& \longrightarrow
& u_1u_2\cdots u_m=0 \ \text{ on }\  \mathbb S_n
\\[6pt]
&& \big\downarrow
\\[6pt]
T_\rho
& \longleftarrow
& u_j=\rho\phi_j
\end{array}
$$

Based on the above discussion, we now outline the main steps of the proof of Theorem \ref{thm:main}.
Suppose that the Toeplitz product is zero but none
of its symbols is identically zero.

\begin{enumerate}
\item[(i)]
Section \ref{sec:operators} establishes the singular-value estimates
needed to extract the factor $T_\rho$.
Starting from the decay of $s_k(T_\rho)$, we control
the differences arising from replacing $T_{\rho\varphi}$
by $T_\varphi T_\rho$ and moving $T_\rho$ past
Toeplitz operators with continuous symbols.
Proposition \ref{thm:product-reduction} consequently shows that a product containing $r$ factors of the form $T_{\rho\varphi}$ can be written as $T_\psi (T_\rho)^r + E$, where $s_k(E) = o(k^{-\frac{r}{n}})\ (k\rightarrow \infty)$.\\

\item[(ii)]
To distinguish this principal term from the remainder,
Proposition \ref{lem:lower-bound} shows that $T_\psi(T_\rho)^r$
cannot have singular values $o(k^{-\frac{r}{n}})$
when $\psi$ is not identically zero on the sphere.
Therefore, the principal term cannot be canceled by a
remainder with the decay obtained in (i).\\

\item[(iii)]
Section \ref{sec:local} makes these estimates applicable to local
boundary factorizations.
Lemma \ref{lem:harmonic} gives $u=\rho\Phi$ near a boundary open set
on which a nonzero harmonic symbol vanishes.
It also ensures that $\Phi$ is nowhere zero on a
smaller boundary open set.
Then we extend these local functions to the closed  unit ball.
Proposition \ref{lem:localization} shows that, after localization,
changes to the symbols away from the chosen set
produce differences whose singular values decay
faster than every negative power.\\

\item[(iv)]
Finally, since the Toeplitz product is zero,
Coburn's theorem  gives that $u_1u_2\cdots u_m=0$ on  $\mathbb S_n$ (see Step 1 in Section \ref{sec:main-proof}). On a suitable boundary open set, the factorizations
in (iii) provide $r \geqslant 1$ factors carrying $\rho$, with the function $\psi$ in the principal term $T_\psi (T_{\rho})^r$ nonzero on a boundary open set.
Then localization and the zero-product assumption
force the principal term
to have singular values $o(k^{-\frac{r}{n}})\ (k\rightarrow \infty)$.
This contradicts (ii).
\end{enumerate}

\section{Preliminaries}

\subsection{Homogeneous decomposition of the Bergman space $L_a^2(\mathbb B_n)$}\label{sec:index-decomposition}\mbox{}

In this subsection, we first introduce the multi-index notation.  Let $\mathbb N_0=\{0,1,2,\cdots\}$. For $\alpha=(\alpha_1,\alpha_2, \cdots,\alpha_n)\in\mathbb N_0^n$ and $z=(z_1, z_2, \cdots,z_n)\in\mathbb C^n$, we write
$$
|\alpha|=\alpha_1+\alpha_2+\cdots+\alpha_n,
\qquad
\alpha!=\alpha_1!\alpha_2!\cdots\alpha_n!,
\qquad
z^\alpha=z_1^{\alpha_1} z_2^{\alpha_2} \cdots  z_n^{\alpha_n},
$$
where  $0!=1$.

With respect to the normalized volume measure $dV$, the monomials satisfy
$$
\|z^\alpha\|_2^2=\frac{n!\alpha!}{(|\alpha|+n)!}.
$$
For each $\alpha\in\mathbb N_0^n$, define the normalized monomial $e_\alpha\colon\mathbb B_n\to\mathbb C$ by
\begin{equation}\label{eq:normalized monomial}
e_\alpha(z)
=\frac{z^\alpha}{\|z^\alpha\|_2}
=\left(\frac{(|\alpha|+n)!}{n!\alpha!}\right)^{\frac{1}{2}}z^\alpha,
\qquad z\in\mathbb B_n.
\end{equation}
The family $\{e_\alpha\}_{\alpha\in\mathbb N_0^n}$ is an orthonormal basis for $L_a^2(\mathbb B_n)$. The above notation and results are taken from  \cite{Coburn1973}.

Since the normalized monomials form an orthonormal basis, every $f\in L_a^2(\mathbb B_n)$ admits the expansion
$$
f
=
\sum_{\alpha\in\mathbb N_0^n}
\langle f,e_\alpha\rangle e_\alpha,
$$
with convergence in $L_a^2(\mathbb B_n)$. Furthermore, we have
$$
\|f\|_2^2
=
\sum_{\alpha\in\mathbb N_0^n}
|\langle f,e_\alpha\rangle|^2.
$$

For each $\ell\in\mathbb N_0$, let
$$
\mathcal H_\ell
=
\operatorname{span}\{z^\alpha:|\alpha|=\ell\}
$$
be the space of homogeneous holomorphic polynomials
of degree $\ell$. Grouping the basis expansion by
total degree, we obtain
$$
f=\sum_{\ell=0}^\infty f_\ell,
\qquad
f_\ell
=
\sum_{|\alpha|=\ell}
\langle f,e_\alpha\rangle e_\alpha
\in \mathcal H_\ell.
$$
The spaces $\mathcal H_\ell\ (\ell \in \mathbb N_0)$ are pairwise orthogonal, and
$$
\left\|f-\sum_{\ell=0}^N f_\ell\right\|_2^2
=
\sum_{|\alpha|>N}
|\langle f,e_\alpha\rangle|^2
\rightarrow0
\ \ \ \  (N\to\infty).
$$
Consequently,
$$
L_a^2(\mathbb B_n)
=
\bigoplus_{\ell=0}^\infty \mathcal H_\ell,
\qquad
\|f\|_2^2
=
\sum_{\ell=0}^\infty\|f_\ell\|_2^2,
$$
where the above direct sum is taken in the Hilbert-space sense. The dimension of $\mathcal H_\ell$ equals the number of multi-indices $\alpha\in\mathbb N_0^n$ satisfying $|\alpha|=\ell$. Thus
$$\dim (\mathcal H_\ell)=\binom{\ell+n-1}{n-1},$$
where $\binom{k}{j}$ stands for the binomial coefficient $\frac{k!}{j!(k-j)!}$.
Similarly, counting the multi-indices with $|\alpha|\leqslant N$ gives
$$
\sum_{\ell=0}^{N}\dim (\mathcal H_\ell)=\binom{N+n}{n}.
$$
For each $N\in\mathbb N_0$, define $Q_N$ on $L_a^2(\mathbb B_n)$ by
\begin{align*}\label{QN}
Q_Nf=\sum_{\ell=0}^{N}f_\ell=\sum_{|\alpha|\leqslant N}\langle f,e_\alpha\rangle e_\alpha=\sum_{|\alpha|\leqslant N}(e_\alpha\otimes e_\alpha) f,
\end{align*}
where $e_\alpha\otimes e_\alpha$ denotes the rank-one operator given by
$(e_\alpha\otimes e_\alpha)f=\langle f, e_\alpha\rangle e_\alpha$.
Then $Q_N$ is the orthogonal projection from $L_a^2(\mathbb B_n)$ onto the space of holomorphic polynomials of total degree at most $N$, namely
$\mathcal H_0\oplus \mathcal H_1\oplus\cdots\oplus\mathcal H_N.$
For each integer $N\geqslant1$, define
$$
\Pi_N f=\sum_{\ell=N}^{2N}f_\ell
=\sum_{N\leqslant|\alpha|\leqslant2N}
\langle f,e_\alpha\rangle e_\alpha.
$$
Then $\Pi_N$ is the orthogonal projection onto
$\mathcal H_N\oplus \mathcal H_{N+1}\oplus\cdots\oplus \mathcal H_{2N}$. Setting
$D_N=\operatorname{rank}(\Pi_N)$,  we have
\begin{equation}\label{eq:band-dimension}
\begin{aligned}
D_N
&=\sum_{\ell=N}^{2N}\dim (\mathcal H_\ell)\\
&=\binom{2N+n}{n}-\binom{N+n-1}{n}\\
&=\frac{2^n-1}{n!}N^n+O(N^{n-1})
\qquad(N\to\infty).
\end{aligned}
\end{equation}

The following lemma allows us to estimate the contribution
of compact operators on the ranges of $\Pi_N$.
\begin{lemma}\label{lem:compact-band}
Let $C$ be a compact operator on $L_a^2(\mathbb B_n)$.
Then
$$
\lim_{N\to\infty}\|C\Pi_N\|=0.
$$
\end{lemma}

\begin{proof}
For $f=\sum\limits_{\ell=0}^{\infty}f_\ell\in L_a^2(\mathbb B_n)$,
orthogonality gives
$$
\|\Pi_Nf\|_2^2
=\sum_{\ell=N}^{2N}\|f_\ell\|_2^2
\leqslant\sum_{\ell=N}^{\infty}\|f_\ell\|_2^2
\rightarrow0
\qquad(N\to\infty).
$$

We first consider a finite-rank operator $F$ on $L_a^2(\mathbb B_n)$.
Write
$$
Ff=\sum_{j=1}^{d}\langle f,g_j\rangle h_j,
\qquad g_j,h_j\in L_a^2(\mathbb B_n).
$$
Since $\Pi_N$ is self-adjoint, we have that
\begin{align*}
\|F\Pi_Nf\|_2
&=\left\|
\sum_{j=1}^{d}\langle f,\Pi_Ng_j\rangle h_j
\right\|_2
\leqslant
\sum_{j=1}^{d}
|\langle f,\Pi_Ng_j\rangle|\,\|h_j\|_2
\leqslant
\sum_{j=1}^{d}\|\Pi_Ng_j\|_2\|h_j\|_2
\end{align*}
for every $\|f\|_2\leqslant1$.
It follows that
$$
\|F\Pi_N\|
\leqslant
\sum_{j=1}^{d}\|\Pi_Ng_j\|_2\|h_j\|_2
\rightarrow0\qquad (N\rightarrow \infty).
$$

Given $\varepsilon>0$, choose a finite-rank operator
$F$ such that $\|C-F\|<\varepsilon$.
Using $\|\Pi_N\|=1$ for all $N$, we obtain
\begin{align*}
\|C\Pi_N\|
\leqslant\|(C-F)\Pi_N\|+\|F\Pi_N\|
\leqslant\|C-F\|+\|F\Pi_N\|
<\varepsilon+\|F\Pi_N\|.
\end{align*}
Since $\varepsilon>0$ is arbitrary, keeping $F$ fixed and letting $N\to\infty$ gives
$$\lim_{N\to\infty}\|C\Pi_N\|=0,$$
as desired.
\end{proof}

\subsection{Singular values of compact operators} \mbox{}

 Let $K$ be a compact operator on a separable Hilbert space $\mathscr H$. Write the polar decomposition of $K$ as $K=U|K|$, where $|K|=(K^*K)^{\frac{1}{2}}$ and $U$ is a partial isometry. Since $|K|$ is positive and compact, there exists a sequence of real numbers $\{\lambda_j\}_{j\geqslant1}$ such that
$$
|K|f=\sum_{j\geqslant 1}\lambda_j\langle f,e_j\rangle e_j, \ \ \ \ f\in \mathscr H,
$$
where $\{e_j\}_{j\geqslant 1}$ is an orthonormal set of eigenvectors corresponding to the positive eigenvalues $\lambda_j$ of $|K|$, listed in nonincreasing order and counted with multiplicity.

The vectors $g_j=Ue_j$ form an orthonormal set, and
\begin{equation}\label{eq:define singular values}
  Kf=\sum_{j\geqslant 1}\lambda_j\langle f,e_j\rangle g_j, \qquad f\in \mathscr H.
\end{equation}
This is the canonical decomposition of $K$, and the numbers $\lambda_j$ are called its nonzero singular values. Let $s_k(K)$ denote the $k$-th singular value of $K$, with multiplicities included and the ordering
$$
s_1(K)\geqslant s_2(K)\geqslant\cdots\geqslant 0.
$$
Taking adjoints in \eqref{eq:define singular values} gives
$$K^*f=\sum_{j\geqslant 1}\lambda_j\langle f,g_j\rangle e_j, \qquad f\in \mathscr H,$$
which is a singular-value decomposition of $K^*$. Therefore $K$ and $K^*$ have the same nonzero singular values, including multiplicities. Appending zeros when the rank is finite, we obtain
\begin{equation}\label{eq:s_k(K^*)=s_k(K)}
s_k(K^*)=s_k(K),\qquad k\geqslant1.
\end{equation}
In particular, if $K$ is positive and compact, then $K^*=K$, and
$$
|K|=(K^*K)^{\frac{1}{2}}=(K^2)^{\frac{1}{2}}=K.
$$
Consequently, the nonzero singular values of $K$ are its positive eigenvalues, listed in nonincreasing order and counted with multiplicity.
For more information on the singular values of  compact operators, one can consult the book \cite{Zhu-book}.

The following lemma expresses the sum of the squared singular values of a compact operator in terms of its action on an orthonormal basis.
\begin{lemma}[{\cite[Theorem~1.22]{Zhu-book}}]
\label{lem:singular-value-square-sum}
Let $K$ be a compact operator on a separable Hilbert space $\mathscr H$,
and let $\{e_j:j\geqslant1\}$ be an orthonormal basis
of $\mathscr H$. Then
\begin{equation*}\label{eq:singular-value-square-sum}
\sum_{k=1}^{\infty}[s_k(K)]^2
=
\sum_{j=1}^{\infty}\|Ke_j\|^2.
\end{equation*}
\end{lemma}

We also need the following lemma to estimate the singular values of some Toeplitz operators on $L_a^2(\mathbb B_n)$.

\begin{lemma}[{\cite[Theorem~1.34]{Zhu-book}}] \label{lem:s_k chara}
Let $\mathscr H_1$ and $\mathscr  H_2$ be separable Hilbert spaces.
  Suppose $K:\mathscr H_1\rightarrow \mathscr H_2$ is compact and $\{s_j(K)\}_{j\geqslant 1}$ is the singular value sequence of $K$. For each $k\geqslant0$, we have
\begin{equation}\label{eq:approximation-number}
s_{k+1}(K)=\inf\{\|K-F\|: F\in \mathcal{F}_k\},
\end{equation}
where $\mathcal{F}_k$ is the set of all bounded linear operators from $\mathscr H_1$ to $\mathscr H_2$ with rank less than or equal to $k$.
\end{lemma}

\begin{remark}\label{rm:compact}
Let $\mathscr H_1$ and $\mathscr H_2$ be separable Hilbert spaces, and let $K:\mathscr H_1\to\mathscr H_2$ be a compact operator. For any $\varepsilon>0$, there exists a finite-rank operator $F$ such that
$$\|K-F\|<\varepsilon.$$
If $\operatorname{rank}(F)< k$, then the above lemma gives
$$
s_k(K)\leqslant \inf\{\|K-G\|: \operatorname{rank}(G)<k\}\leqslant \|K-F\|<\varepsilon.
$$
Hence $s_k(K)\to0$ as $k\to\infty$.
\end{remark}

From now on, when $x$ is a real number, we shall use $\lfloor x \rfloor$ to denote the greatest integer not exceeding $x$.

\begin{lemma}\label{lem:singular-value inequalities}
Assume that each  $\mathscr H_j$ is a separable  Hilbert space, where $j=0, 1, 2, 3$.
Suppose that $A: \mathscr H_0\rightarrow \mathscr H_1$  and $B: \mathscr H_2\rightarrow \mathscr H_3$ are both  bounded and linear, and $C: \mathscr H_1\rightarrow \mathscr H_2$, $K: \mathscr H_1\rightarrow \mathscr H_2$, $K': \mathscr H_1\rightarrow \mathscr H_2$ and
 $K'': \mathscr H_2\rightarrow \mathscr H_3$ are all compact. Then the following three assertions hold:
\begin{enumerate}
\item[\emph{(a)}]
For every integer $k\geqslant1$,
\begin{equation}\label{eq:sv-sum}
s_{2k-1}(C+K')
\leqslant s_k(C)+s_k(K')
\end{equation}
and
\begin{equation}\label{eq:sv-product}
s_{2k-1}(K''K)
\leqslant s_k(K'')s_k(K).
\end{equation}

\item[\emph{(b)}]
For every integer $k\geqslant 1$,
\begin{equation}\label{eq:sv-bounded}
s_k(BKA)
\leqslant \|A\|\,\|B\|s_k(K).
\end{equation}

\item[\emph{(c)}]
If
$$
s_k(K)=O(k^{-\gamma})
\qquad (k\to\infty)
$$
for some $\gamma>0$, then
\begin{equation}\label{eq:sv-compact-product}
s_k(K''K)=o(k^{-\gamma})
\qquad (k\to\infty).
\end{equation}
\end{enumerate}
\end{lemma}

\begin{proof}
Inequalities \eqref{eq:sv-sum} and
\eqref{eq:sv-product} follow from
 \cite[Corollary 1.35]{Zhu-book} immediately. We next prove (b).
For every operator $F$ with $\operatorname{rank}(F)<k$, we have
$$\operatorname{rank}(BFA)\leqslant \operatorname{rank}(F)<k.$$
Applying Lemma \ref{lem:s_k chara} to the compact operator $BKA$, we obtain
$$
\begin{aligned}
s_k(BKA)
\leqslant \|BKA-BFA\|
= \|B(K-F)A\|
\leqslant \|A\|\,\|B\|\,\|K-F\|.
\end{aligned}
$$
Taking the infimum over all such $F$ gives (b).

It remains to prove (c). By \eqref{eq:sv-product},
$$
s_{2k-1}(K''K)\leqslant s_k(K'')s_k(K).
$$
Since $K''$ is compact, $s_k(K'')\to 0$ as $k\rightarrow \infty$.
The assumption on $K$ implies that
$k^\gamma s_k(K)$ is bounded for sufficiently
large $k$. Hence
$$
k^\gamma s_{2k-1}(K''K)
\leqslant
s_k(K'')\bigl[k^\gamma s_k(K)\bigr]
\rightarrow0
\qquad (k\to\infty).
$$

For an arbitrary integer $l\geqslant1$, let
$$
k=\Big\lfloor \frac{l+1}{2}\Big\rfloor.
$$
Then we have that $2k-1\leqslant l\leqslant2k$ and
$k\to\infty$ as $l\to\infty$.
By the monotonicity of singular values,
$$
\begin{aligned}
l^\gamma s_l(K''K)
\leqslant (2k)^\gamma s_{2k-1}(K''K)
\leqslant
2^\gamma s_k(K'')\bigl[k^\gamma s_k(K)\bigr]
\rightarrow0
\qquad (l\to\infty).
\end{aligned}
$$
This proves \eqref{eq:sv-compact-product}, to complete the proof of Lemma \ref{lem:singular-value inequalities}.
\end{proof}

Let $K_1, K_2,\cdots, K_q$ be compact operators, where $q$ is a fixed positive integer.  Suppose that
$$s_k(K_j)=O(k^{-\gamma_j})\qquad (k\to\infty)$$
for some  $\gamma_j>0$ and $1\leqslant j\leqslant q$.
Let $\gamma=\min\limits_{1\leqslant j\leqslant q}\gamma_j.$
For any integer $l\geqslant1$, let
$$k=\Big\lfloor\frac{l-1}{2^{q-1}}\Big\rfloor+1.$$
Then $2^{q-1}(k-1)+1\leqslant l\leqslant 2^{q-1}k$.  Repeated application of (\ref{eq:sv-sum}) in Lemma~\ref{lem:singular-value inequalities} gives
$$s_{2^{q-1}(k-1)+1}(K_1+K_2+\cdots+K_q)\leqslant\sum_{j=1}^{q} s_k(K_j).$$
As $l\leqslant 2^{q-1}k$, we have
$$
\begin{aligned}
l^\gamma s_l(K_1+K_2+\cdots+K_q)
&\leqslant (2^{q-1}k)^\gamma s_{2^{q-1}(k-1)+1}(K_1+K_2+\cdots+K_q)\\
&\leqslant 2^{\gamma(q-1)}\sum_{j=1}^q k^\gamma s_k(K_j)\\
&= 2^{\gamma(q-1)}\sum_{j=1}^q k^{\gamma-\gamma_j}\big[k^{\gamma_j} s_k(K_j)\big].
\end{aligned}
$$
Therefore,
\begin{align}\label{finitesum}
s_l(K_1+K_2+\cdots+K_q)=O(l^{-\gamma})\qquad (l\to\infty).
\end{align}
From the above discussion, we  further obtain  that
\begin{equation}\label{eq:s_k-finite-sum}
  s_l(K_1+K_2+\cdots+K_q)=o(l^{-\gamma})\qquad (l\to\infty)
\end{equation}
if
$$s_k(K_j)=o(k^{-\gamma_j})\qquad (k\to\infty)$$
for every $j$.

We next consider the product $$B_1 B_2\cdots B_s K_1 K_2\cdots K_t,$$
where $B_1, B_2, \cdots,B_s$ are bounded, $K_1, K_2\cdots, K_t$ are compact, and the domain and range of adjacent factors are compatible. Assume that
$$
s_k(K_j)=O(k^{-\gamma_j}) \qquad(k\to\infty)$$ and $$ \gamma_j>0,\quad 1\leqslant j\leqslant t.
$$
Applying \eqref{eq:sv-bounded} and then repeatedly applying (\ref{eq:sv-product}), we obtain
$$
s_{2^{t-1}(k-1)+1}\bigl(B_1 B_2\cdots B_sK_1K_2\cdots K_t\bigr)\leqslant\bigg(\prod_{j=1}^{s}\|B_j\|\bigg)\prod_{j=1}^{t}s_k(K_j).
$$
Similarly, for any integer $l\geqslant1$, let
$$k=\Big\lfloor\frac{l-1}{2^{t-1}}\Big\rfloor+1.$$
Then $2^{t-1}(k-1)+1\leqslant l\leqslant 2^{t-1}k$ and
$$
l^{\gamma_1+\gamma_2+\cdots+\gamma_t}s_l\bigl(B_1B_2\cdots B_sK_1K_2\cdots K_t\bigr)\leqslant 2^{(t-1)(\gamma_1+\gamma_2+\cdots+\gamma_t)}
\bigg(\prod_{j=1}^{s}\|B_j\|\bigg)\prod_{j=1}^{t}\bigl[k^{\gamma_j}s_k(K_j)\bigr].
$$
Since each factor $k^{\gamma_j}s_k(K_j)$ is bounded for sufficiently large $k$, we have
\begin{equation}\label{eq:product-s_k-equ}
s_l\bigl(B_1B_2\cdots B_sK_1K_2\cdots K_t\bigr)=O\bigl(l^{-(\gamma_1+\gamma_2+\cdots+\gamma_t)}\bigr)\qquad(l\to\infty).
\end{equation}
If at least one $K_j$ satisfies $s_k(K_j)=o(k^{-\gamma_j})$,
then the corresponding factor $k^{\gamma_j}s_k(K_j)$ tends to zero. Consequently,
\begin{equation}\label{eq:s_k product estimate-small}
  s_l\bigl(B_1B_2\cdots B_sK_1K_2\cdots K_t\bigr)=o\bigl(l^{-(\gamma_1+\gamma_2+\cdots+\gamma_t)}\bigr)\qquad(l\to\infty).
\end{equation}

\subsection{The Laplacian in polar coordinates}\mbox{}

In the sequel, we identify $\mathbb C^n$ with $\mathbb R^{2n}$ and write $z_j=x_j+\mathrm{i}y_j$. Define the Euclidean Laplacian by
$$\Delta u=\sum_{j=1}^{n}\left(\frac{\partial^2u}{\partial x_j^2}+\frac{\partial^2u}{\partial y_j^2}\right).$$
A function $u\in C^2(\mathbb B_n)$ (i.e., $u$  is twice continuously differentiable on $\mathbb B_n$) is harmonic if and only if $\Delta u=0$ on $\mathbb B_n$.

Note that each nonzero point in $\mathbb B_n$ can be written as $t\zeta$, where $0<t<1$ and $\zeta\in\mathbb S_n$. For fixed $t$, the values $u(t\zeta)$, with $\zeta\in\mathbb S_n$, define a function on the unit sphere.
For fixed $\zeta$, varying $t$ moves the point $t\zeta$ along the ray from the origin in direction $\zeta$. Accordingly,
$
\frac{\partial}{\partial t}u(t\zeta)$ and
$\frac{\partial^2}{\partial t^2}u(t\zeta)
$
are called the first and second radial derivatives of $u$, respectively.

The unit sphere $\mathbb S_n$ carries the metric induced by the Euclidean metric. The Laplace--Beltrami operator is defined by
$$
\Delta_{\mathbb S_n}u
=
\operatorname{div}_{\mathbb S_n}
\bigl(\nabla_{\mathbb S_n}u\bigr),
$$
where  \(\nabla_{\mathbb S_n}\) is the gradient on the sphere and \(\operatorname{div}_{\mathbb S_n}\) is the divergence of a tangent vector field on the sphere. Equivalently, for $u\in C^2(\mathbb S_n)$,
extend $u$ to $\mathbb C^n\backslash\{0\}$ by
$z\mapsto u\big(\frac{z}{|z|}\big)$.
Then the sphere Laplacian is given by
\begin{equation*}\label{eq:sphere-laplacian-definition}
\Delta_{\mathbb S_n}u(\zeta)
=
\left.
\Delta\left[u\left(\frac{z}{|z|}\right)\right]
\right|_{z=\zeta},
\qquad \zeta\in\mathbb S_n.
\end{equation*}

We recall the polar-coordinate formula for the
Laplacian, stated below in our notation.

\begin{lemma}[{\cite[Lemma~1.4.1]{DaiXu}}]
\label{lem:polar-laplacian}
Let $u\in C^2(\mathbb B_n)$.
For $0<t<1$ and $\zeta\in\mathbb S_n$,
\begin{equation*}
\Delta u(t\zeta)
=
\frac{\partial^2}{\partial t^2}u(t\zeta)
+
\frac{2n-1}{t}
\frac{\partial}{\partial t}u(t\zeta)
+
\frac{1}{t^2}
\Delta_{\mathbb S_n}u(t\zeta).
\end{equation*}
In particular, the above equation  becomes
\begin{equation*}\label{eq:prelim-polar-harmonic}
\frac{\partial^2}{\partial t^2}u(t\zeta)
+
\frac{2n-1}{t}
\frac{\partial}{\partial t}u(t\zeta)
+
\frac{1}{t^2}
\Delta_{\mathbb S_n}u(t\zeta)
=
0
\end{equation*}
when  $u$ is harmonic on $\mathbb B_n$.
\end{lemma}

\section{Singular-value estimates for Toeplitz products}\label{sec:operators}
In the following, we denote
$$ \rho(z)=1-|z|^2.$$
In this section, we first compute the action of the Toeplitz operator $T_\rho$ on the homogeneous spaces and estimate the singular values of $T_{\rho \varphi}$ with $\varphi\in L^\infty(\mathbb B_n)$; then we use these estimates to reduce a finite Toeplitz product to a term of the form $T_{\psi}(T_\rho)^r$; finally, we show that this term cannot have singular values of order $o(k^{-\frac{r}{n}})$ when the boundary values of $\psi$ are not identically zero.

Let us begin by studying the Toeplitz operator $T_{\rho}$.

\begin{lemma}\label{lem:radial}
Let $\rho(z)=1-|z|^2$. The Toeplitz operator $T_\rho$ is positive, injective and compact on $L_a^2(\mathbb B_n)$. Moreover, for $\ell\geqslant 0$ we have
\begin{equation}\label{eq:R-spectrum}
T_\rho f_\ell=\frac1{\ell+n+1}f_\ell,\qquad f_\ell\in\Hh_\ell.
\end{equation}
\end{lemma}

\begin{proof}
We first compute the action of $T_\rho$ on each homogeneous space. If $f_\ell\in\Hh_\ell$, then
$$f_\ell(z)=\sum_{|\alpha|=\ell}c_\alpha z^\alpha.$$
Regarding $f_\ell$ as the homogeneous polynomial defined on $\mathbb C^n$,  we have
\begin{align}\label{f_l(tz)=t^lf(z)}
  \begin{split}
    f_\ell(t\zeta) &=\sum_{|\alpha|=\ell}c_\alpha(t\zeta_1)^{\alpha_1}(t\zeta_2)^{\alpha_2}\cdots(t\zeta_n)^{\alpha_n}\\
    &=\sum_{|\alpha|=\ell}c_\alpha t^{\alpha_1+\alpha_2+\cdots+\alpha_n}\zeta^\alpha\\
    &=t^\ell f_\ell(\zeta)
  \end{split}
\end{align}
 for any $t\in \mathbb{C}$ and $\zeta=(\zeta_1, \zeta_2, \cdots,\zeta_n)\in \mathbb C^n$.

Let $g_j\in\mathcal H_j$. Since $P$ is the orthogonal projection onto $L_a^2(\mathbb B_n)$, we have
\begin{align*}\label{inner product,f_l,g_j}
  \begin{split}
     \langle T_\rho f_\ell,g_j\rangle &=\langle P(\rho f_\ell),g_j\rangle
     =\langle \rho f_\ell ,g_j\rangle
     =\int_{\B}(1-|z|^2)f_\ell(z)\overline{g_j(z)}d V(z).
  \end{split}
\end{align*}
For any $\theta\in\mathbb R$, it follows from \eqref{f_l(tz)=t^lf(z)} that
\[
f_\ell(e^{\mathrm i\theta}\xi)=e^{\mathrm i\ell\theta}f_\ell(\xi),
\qquad
\overline{g_j(e^{\mathrm i\theta}\xi)}
=e^{-\mathrm ij\theta}\overline{g_j(\xi)}.
\]
Consequently, we have
\begin{align*}
\langle T_\rho f_\ell,g_j\rangle
&=\int_{\mathbb B_n}(1-|\xi|^2)f_\ell(e^{\mathrm i \theta}\xi)\overline{g_j(e^{\mathrm i \theta}\xi)}\,dV(\xi)\\
&=e^{\mathrm i (\ell-j)\theta}\int_{\mathbb B_n}(1-|\xi|^2)f_\ell(\xi)\overline{g_j(\xi)}\,dV(\xi)\\
&=e^{\mathrm i(\ell-j)\theta}\langle T_\rho f_\ell,g_j\rangle.
\end{align*}

If $\ell\neq j$, then we take $\theta=\frac{\pi}{|\ell-j|}$ to obtain that
$e^{\mathrm i(\ell-j)\theta}=-1$, and hence
$$\langle T_\rho f_\ell,g_j\rangle=-\langle T_\rho f_\ell,g_j\rangle.$$
It follows that
$$\langle T_\rho f_\ell,g_j\rangle=0, \qquad j\neq\ell.$$
In the case  of $\ell=j$, we have by the homogeneity that
$$f_\ell(t\zeta)\overline{g_\ell(t\zeta)}=t^{2\ell}f_\ell(\zeta)\overline{g_\ell(\zeta)}$$
for
$0<t<1$ and $\zeta\in\mathbb S_n$.
Using  \eqref{polar}, we obtain
\begin{align*}
\langle T_\rho f_\ell,g_\ell\rangle &=\int_{\mathbb B_n}(1-|z|^2) f_\ell(z)\overline{g_\ell(z)}\, dV(z)\\
&=2n\left(\int_0^1(1-t^2)t^{2\ell+2n-1}\, dt\right)\left(\int_{\mathbb S_n}f_\ell(\zeta)\overline{g_\ell(\zeta)}\, d\sigma(\zeta)\right)
\end{align*}
and
$$\langle f_\ell,g_\ell\rangle=2n\left(\int_0^1t^{2\ell+2n-1}\,dt\right)\left(\int_{\mathbb S_n}
f_\ell(\zeta)\overline{g_\ell(\zeta)}\,d\sigma(\zeta)\right).
$$
Since
$$\int_0^1t^{2\ell+2n-1}\,dt=\frac{1}{2(\ell+n)}$$
and
\begin{align*}
\int_0^1(1-t^2)t^{2\ell+2n-1}\,dt=
\frac{1}{2(\ell+n)(\ell+n+1)},
\end{align*}
 we obtain that
$$\langle T_\rho f_\ell,g_\ell\rangle=\frac{1}{\ell+n+1}\langle f_\ell,g_\ell\rangle$$
and
$$\langle T_\rho f_\ell,g_j\rangle
=\begin{cases}
\frac{1}{\ell+n+1} \langle f_\ell,g_j\rangle,&j=\ell,\vspace{2mm}\\
0,&j\ne\ell.
\end{cases}
$$
Since \(f_\ell\) is orthogonal to \(g_j\) for \(\ell \neq j\) and the spaces of homogeneous polynomials span a dense subspace, this proves \eqref{eq:R-spectrum}.

For every $f\in L_a^2(\mathbb B_n)$,
$$
\langle T_\rho f,f\rangle
=\int_{\mathbb B_n}\rho(z)|f(z)|^2\,dV(z)\geqslant0.
$$
Thus $T_\rho$ is positive. If $T_\rho f=0$, then this integral is zero. Since $\rho>0$ on $\mathbb B_n$, it follows that $f=0$. Hence $T_\rho$ is injective on $L_a^2(\mathbb B_n)$. Finally, $T_\rho$ is diagonal with respect to the orthonormal
basis $\{e_\alpha\}_{\alpha\in\mathbb N_0^n}$  of this space.  Each homogeneous space is finite dimensional, and
the diagonal entries tend to zero as $|\alpha|\rightarrow \infty$. Thus, we have that
$T_\rho$ is compact on $L_a^2(\mathbb B_n)$.
\end{proof}

The next proposition gives the asymptotic estimate for the singular values of the Toeplitz operators \(T_\rho\)  and $T_{\rho \varphi}$ with $\varphi \in L^\infty(\mathbb B_n)$.

\begin{proposition}\label{lem:s_k(T_rho)-eqvi and s_k(T_rho)-leq}
Let $\rho(z)=1-|z|^2$. We have
\begin{equation*}\label{eq:s_k(T_rho)-leq}
  s_k(T_\rho)=O( k^{-\frac 1 n})\qquad(k\to\infty).
\end{equation*}
Moreover, for every $\varphi \in L^\infty(\mathbb B_n)$, the Toeplitz operator $T_{\rho \varphi}$ is compact and
  \begin{align*}\label{eq:s_k(T_a rho)}
s_k(T_{\rho \varphi})=O(k^{-\frac{1}{n}})\qquad(k\to\infty).
\end{align*}
\end{proposition}

\begin{proof}
By Lemma \ref{lem:radial}, $T_\rho$ is positive and compact. Its singular values are therefore its positive eigenvalues. Then \eqref{eq:R-spectrum} shows that the eigenvalue
$$\frac1{\ell+n+1}$$
has multiplicity $\dim (\mathcal H_\ell)$.
We first count the eigenvalues contributed by the homogeneous spaces. Using the identity
$$
\binom{j+n-1}{n-1}=\binom{j+n}{n}-\binom{j+n-1}{n}
$$
and summing over $0\leqslant j\leqslant\ell$, we obtain
$$
\begin{aligned}
\sum_{j=0}^{\ell}\dim (\mathcal H_j)&=\sum_{j=0}^{\ell}
\left[\binom{j+n}{n}-\binom{j+n-1}{n}\right]=\binom{\ell+n}{n},
\end{aligned}
$$
where we adopt the convention that \(\binom{n-1}{n}=0\).  Similarly,
$$
\sum_{j=0}^{\ell-1}\dim (\mathcal H_j)=\binom{\ell+n-1}{n}.
$$
Thus, for each $k\geqslant1$, there is a unique integer $\ell_k\geqslant0$ such that
\begin{equation}\label{eq:ineq of k}
\binom{\ell_k+n-1}{n}<k\leqslant\binom{\ell_k+n}{n}<\frac{(\ell_k+n+1)^n}{n!}.
\end{equation}
The first inequality says that the eigenvalues contributed by degrees below $\ell_k$ occupy fewer than $k$ positions, and  the second inequality says that including degree $\ell_k$ reaches the $k$-th position.
Therefore, by (\ref{eq:ineq of k}) we have  that $\ell_k\to\infty$ as $k\rightarrow \infty$, and
\begin{equation}\label{eq:skT}
  s_k(T_\rho)=\frac{1}{\ell_k+n+1}.
\end{equation}
For $\ell_k\geqslant1$, expanding the binomial coefficients and dividing by $\ell_k^n$, we obtain
$$
\frac{1}{n!}
\prod_{j=0}^{n-1}\left(1+\frac{j}{\ell_k}\right)
<\frac{k}{\ell_k^n}\leqslant\frac{1}{n!}\prod_{j=1}^{n}\left(1+\frac{j}{\ell_k}\right),
$$
which gives
$$
\lim_{k\to\infty}\frac{k}{\ell_k^n}=\lim_{k\to\infty}\frac{1}{n!}\prod_{j=0}^{n-1}\left(1+\frac{j}{\ell_k}\right)=\lim_{k\to\infty}\frac{1}{n!}\prod_{j=1}^{n}\left(1+\frac{j}{\ell_k}\right)=\frac{1}{n!}.
$$
Using  \eqref{eq:skT}, we conclude that
$$
\begin{aligned}
\lim_{k\to\infty}k^{\frac 1 n}s_k(T_\rho)&=\lim_{k\to\infty}\frac{k^{\frac 1 n}}{\ell_k+n+1}
=\lim_{k\to\infty}\left(\frac{k}{\ell_k^n}\right)^{\frac 1 n}\frac{\ell_k}{\ell_k+n+1}=(n!)^{-\frac 1 n}.
\end{aligned}
$$

We next estimate $s_k(T_{\rho \varphi})$. For $f_\ell\in\mathcal H_\ell$, the polar integration formula gives
$$
\|f_\ell\|_2^2=2n\int_0^1 t^{2\ell+2n-1}\,dt\int_{\mathbb S_n}|f_\ell(\zeta)|^2\,d\sigma(\zeta)
=\frac{n}{\ell+n}\int_{\mathbb S_n}|f_\ell(\zeta)|^2\,d\sigma(\zeta)
$$
and
$$
\begin{aligned}
\|\rho f_\ell\|_2^2
&=2n\int_0^1(1-t^2)^2t^{2\ell+2n-1}\,dt\int_{\mathbb S_n}|f_\ell(\zeta)|^2\,d\sigma(\zeta)\\
&=2n\left[\frac{1}{2(\ell+n)}-\frac{2}{2(\ell+n+1)}
+\frac{1}{2(\ell+n+2)}\right]\int_{\mathbb S_n}|f_\ell(\zeta)|^2\,d\sigma(\zeta)\\
&=\frac{2n}{(\ell+n)(\ell+n+1)(\ell+n+2)}\int_{\mathbb S_n}|f_\ell(\zeta)|^2\,d\sigma(\zeta)\\
&=\frac{2}{(\ell+n+1)(\ell+n+2)}\|f_\ell\|_2^2.
\end{aligned}
$$
Combining this  with \eqref{eq:R-spectrum} yields
\begin{equation}\label{eq:rho-homogeneous-norm-bound}
\|\rho f_\ell\|_2^2
\leqslant
\frac{2}{(\ell+n+1)^2}\|f_\ell\|_2^2
=
2\|T_\rho f_\ell\|_2^2.
\end{equation}

To extend \eqref{eq:rho-homogeneous-norm-bound} to an arbitrary function in $L_a^2(\mathbb B_n)$, we first verify that the functions $\rho f_\ell~(\ell \in\mathbb N_0)$ are pairwise orthogonal in $L^2(\mathbb B_n)$. If $\ell\ne j$, the substitution $z=e^{\mathrm i\theta}w$ with $\theta\in \mathbb R$ gives
\[
\begin{aligned}
\langle \rho f_\ell,\rho f_j\rangle
&=
\int_{\mathbb B_n}
\rho(z)^2 f_\ell(z)\overline{f_j(z)}\,dV(z)\\
&=
e^{\mathrm i(\ell-j)\theta}
\int_{\mathbb B_n}
\rho(w)^2 f_\ell(w)\overline{f_j(w)}\,dV(w).
\end{aligned}
\]
Choosing $\theta$ so that $e^{\mathrm i(\ell-j)\theta}=-1$ shows that this inner product is zero, as desired.

Letting
$$f=\sum_{\ell=0}^\infty f_\ell\in L_a^2(\mathbb B_n),$$
it follows from \eqref{eq:rho-homogeneous-norm-bound} that
\begin{equation}\label{eq:radial calculation}
\begin{aligned}
\|\rho f\|_2^2
=\sum_{\ell=0}^\infty\|\rho f_\ell\|_2^2
\leqslant2\sum_{\ell=0}^\infty\frac{\|f_\ell\|_2^2}{(\ell+n+1)^2}
=2\|T_\rho f\|_2^2.
\end{aligned}
\end{equation}
Consequently,
\begin{equation}\label{eq:rho-a-norm-domination}
\begin{aligned}
\|T_{\rho \varphi}f\|_2
\leqslant\|\rho \varphi f\|_2\leqslant\|\varphi \|_\infty\|\rho f\|_2
\leqslant\sqrt{2}\,\|\varphi\|_\infty\|T_\rho f\|_2
\end{aligned}
\end{equation}
for every $\varphi \in L^\infty(\mathbb B_n)$ and $f\in L_a^2(\mathbb B_n)$.

Next, we use this estimate to construct finite-rank
approximations to $T_{\rho \varphi}$. For $\ell\geqslant0$, define
\[
F_\ell f
=
T_{\rho \varphi}
\bigg(\sum_{j=0}^{\ell-1}f_j\bigg),
\qquad
f=\sum_{j=0}^\infty f_j,
\]
where $F_0=0$. The range of $F_\ell$ is contained in
$$
T_{\rho \varphi}(\mathcal H_0\oplus \mathcal H_1\oplus\cdots\oplus\mathcal H_{\ell-1}).
$$
The dimension of this image cannot exceed that of its domain. Hence, we have
\begin{equation}\label{eq:rank-k}
\operatorname{rank}(F_\ell)
\leqslant
\sum_{j=0}^{\ell-1}\dim (\mathcal H_j)
=
\binom{\ell+n-1}{n}.
\end{equation}
Applying \eqref{eq:rho-a-norm-domination}, we obtain
\begin{align*}
\begin{split}
\|(T_{\rho \varphi}-F_\ell)f\|_2^2
&=\bigg\|T_{\rho \varphi}\bigg(\sum_{j=\ell}^\infty f_j\bigg)\bigg\|_2^2
\leqslant2\|\varphi\|_\infty^2\bigg\|T_\rho\bigg(\sum_{j=\ell}^\infty f_j\bigg)\bigg\|_2^2\\
&=2\|\varphi\|_\infty^2\sum_{j=\ell}^\infty\frac{\|f_j\|_2^2}{(j+n+1)^2}
\leqslant\frac{2\|\varphi\|_\infty^2}{(\ell+n+1)^2}\|f\|_2^2.
\end{split}
\end{align*}
Consequently,
\begin{equation*}\label{eq:rho-a-finite-rank-error}
\|T_{\rho \varphi}-F_\ell\|
\leqslant
\frac{\sqrt{2}\,\|\varphi\|_\infty}{\ell+n+1}.
\end{equation*}
The right-hand side tends to zero as $\ell\rightarrow \infty$. Hence $T_{\rho \varphi}$ is compact.

Given $k\geqslant1$, choose $\ell_k$ from \eqref{eq:ineq of k}. Then  \eqref{eq:rank-k} gives $\operatorname{rank}(F_{\ell_k})<k$. By (\ref{eq:skT}) and  Lemma \ref{lem:s_k chara}, we have
$$
\begin{aligned}
s_k(T_{\rho \varphi})&\leqslant\|T_{\rho \varphi}-F_{\ell_k}\|\leqslant\frac{\sqrt{2}\,\|\varphi\|_\infty}{\ell_k+n+1}
=\sqrt{2}\,\|\varphi\|_\infty s_k(T_\rho)\leqslant\sqrt{2}\,(n!)^{-\frac1 n}\|\varphi\|_\infty k^{-\frac 1 n},
\end{aligned}
$$
as desired. This completes the proof.
\end{proof}

In the following, we turn to the estimation of  the singular values of the semi-commutator \(T_{\rho \varphi}-T_{\varphi}T_\rho\) and the commutator \(T_\rho T_{\varphi}-T_{\varphi}T_\rho\).

\begin{lemma}\label{lem:compactness-radial-errors}
 Let $\varphi \in C(\overline{\mathbb B}_n)$. Then
$$
s_k(T_{\rho \varphi}-T_{\varphi}T_\rho)=o(k^{-\frac{1}{n}})\qquad (k\to\infty),
$$
and
$$
s_k(T_\rho T_{\varphi}-T_{\varphi}T_\rho)=o(k^{-\frac{1}{n}})\qquad (k\to\infty).
$$
\end{lemma}

\begin{proof}
We first estimate the singular values of the Hankel operator $H_\rho$. By \eqref{eq:radial calculation}, we have
$$
\|H_\rho f\|_2=\|(I-P)(\rho f)\|_2\leqslant\|\rho f\|_2\leqslant\sqrt{2}\,\|T_\rho f\|_2, \ \ \ \ f\in L_a^2(\mathbb B_n).
$$
Applying this inequality to $(I-Q_N)f$, we obtain
$$
\begin{aligned}
\|H_\rho(I-Q_N)f\|_2^2
\leqslant 2\,\|T_\rho(I-Q_N)f\|_2^2
=2\sum_{\ell=N+1}^{\infty}\frac{\|f_\ell\|_2^2}{(\ell+n+1)^2}
\leqslant\frac{2}{(N+n+2)^2}\|f\|_2^2.
\end{aligned}
$$
Hence
$$
\|H_\rho-H_\rho Q_N\|\leqslant\frac{\sqrt{2}}{N+n+2}
$$
and
$$
\operatorname{rank}(H_\rho Q_N)\leqslant\operatorname{rank}(Q_N)
=\sum_{j=0}^{N}\dim (\mathcal H_j)
=\binom{N+n}{n}.
$$

For each sufficiently large positive integer $k$, choose the unique integer $\ell_k\geqslant 1$ such that
$$
\binom{\ell_k+n-1}{n}<k\leqslant\binom{\ell_k+n}{n}.
$$
Taking $N=\ell_k-1$, we obtain
$$
\operatorname{rank}(H_\rho Q_{\ell_k-1})<k.
$$
The finite-rank characterization of singular values in Lemma \ref{lem:s_k chara} therefore implies
$$
\begin{aligned}
s_k(H_\rho)
\leqslant\|H_\rho-H_\rho Q_{\ell_k-1}\|
\leqslant\frac{\sqrt{2}}{\ell_k+n+1}
=\sqrt{2}\,s_k(T_\rho).
\end{aligned}
$$
By Proposition \ref{lem:s_k(T_rho)-eqvi and s_k(T_rho)-leq} and  \eqref{eq:s_k(K^*)=s_k(K)}, we have
$$
s_k(H_\rho^*)=s_k(H_\rho)=O(k^{-\frac{1}{n}}).
$$

Since $\rho$ is real-valued, the  identity  in \eqref{eq:Hankel identity} yields
$$
T_{\rho \varphi}-T_{\varphi}T_\rho=H_{\overline \varphi}^{\,*}H_\rho
\qquad\text{and}\qquad
T_{\rho \varphi}-T_\rho T_{\varphi}=H_\rho^*H_{\varphi}.
$$
Using the compactness of $H_{\varphi}$ and $H_{\overline \varphi}$, we have by Remark \ref{rm:compact} that
$$
s_k(H_{\varphi})\rightarrow 0
\qquad\text{and}\qquad
s_k(H_{\overline \varphi}^{\,*})\rightarrow 0
$$
as $k\rightarrow \infty$.
Using Lemma \ref{lem:singular-value inequalities}, we obtain
$$
\begin{aligned}
s_{2k-1}(T_{\rho \varphi}-T_{\varphi}T_\rho)
=
s_{2k-1}(H_{\overline \varphi}^{\,*}H_\rho)
\leqslant
s_k(H_{\overline \varphi}^{\,*})s_k(H_\rho)
=o(k^{-\frac{1}{n}}).
\end{aligned}
$$

Similarly, using
$s_k(H_\rho^*)=s_k(H_\rho)$, we have
$$
\begin{aligned}
s_{2k-1}(T_{\rho \varphi}-T_\rho T_{\varphi})
=
s_{2k-1}(H_\rho^*H_\varphi)
\leqslant
s_k(H_\rho^*)s_k(H_\varphi)
=
o(k^{-\frac{1}{n}}).
\end{aligned}
$$
For each $l\geqslant1$, let $k=\lfloor\frac{l+1}{2}\rfloor$. Then we have
$
2k-1\leqslant l\leqslant2k.
$
It follows that
$$
l^{\frac{1}{n}}s_l(T_{\rho \varphi}-T_{\varphi}T_\rho)
\leqslant(2k)^{\frac{1}{n}}s_{2k-1}(T_{\rho \varphi}-T_{\varphi}T_\rho)\rightarrow 0\qquad (k\rightarrow \infty)
$$
and
$$
l^{\frac{1}{n}}s_l(T_{\rho \varphi}-T_\rho T_{\varphi})
\leqslant(2k)^{\frac{1}{n}}s_{2k-1}(T_{\rho \varphi}-T_\rho T_{\varphi})\rightarrow 0\qquad (k\rightarrow \infty).
$$
Thus, we get that
$$
s_k(T_{\rho \varphi}-T_{\varphi}T_\rho)=o(k^{-\frac{1}{n}}) \ \ \ \ \text{and}\ \ \ \
s_k(T_{\rho \varphi}-T_\rho T_{\varphi})=o(k^{-\frac{1}{n}}) \ \ \ \ (k\rightarrow \infty).
$$

Finally, using
$$
T_\rho T_{\varphi}-T_{\varphi}T_\rho
=
(T_{\rho \varphi}-T_{\varphi}T_\rho)
-
(T_{\rho \varphi}-T_\rho T_{\varphi}),
$$
we have by Lemma \ref{lem:singular-value inequalities} that
$$
s_{2k-1}(T_\rho T_{\varphi}-T_{\varphi}T_\rho)\leqslant s_k(T_{\rho \varphi}-T_{\varphi}T_\rho)+s_k(T_{\rho \varphi}-T_\rho T_{\varphi})=o(k^{-\frac{1}{n}}).
$$
Then using the same argument as the one used above,   we arrive at
$$
s_k(T_\rho T_{\varphi}-T_{\varphi}T_\rho)
=
o(k^{-\frac{1}{n}})  \ \ \ \ (k\rightarrow \infty),
$$
as required.
\end{proof}

Next, we  apply Lemma \ref{lem:compactness-radial-errors} to a finite product of Toeplitz operators.

\begin{proposition}\label{thm:product-reduction}
Let $\varphi_1, \varphi_2, \cdots, \varphi_m\in C(\overline{\mathbb B}_n)$.  Choose  $J\subset\{1, 2,\cdots,m\}$,  let $r=|J|$ and define
$$
v_j=
\begin{cases}
\rho \varphi_j, & j\in J,\\
\varphi_j, & j\notin J,
\end{cases}
$$
where $|J|$ denotes the cardinality of the set $J$.
Then
$$
T_{v_1}T_{v_2}\cdots T_{v_m}=T_{\varphi_1\varphi_2\cdots \varphi_m}(T_\rho)^r+E,
$$
where $E$ is compact and
$$
s_k(E)=o(k^{-\frac{r}{n}})\qquad (k\to\infty).
$$
Here $(T_\rho)^0=I$, the identity operator on $L_a^2(\mathbb B_n)$.
\end{proposition}

\begin{proof}
We prove the assertion by induction on the number $m$ of Toeplitz factors. For this purpose, we first establish two estimates which will be used in the induction.

By Proposition \ref{lem:s_k(T_rho)-eqvi and s_k(T_rho)-leq} and  \eqref{eq:product-s_k-equ}, we have
\begin{equation}\label{eq:s_k111111}
  s_k[(T_\rho)^p]=O(k^{-\frac{p}{n}})\ \ \ \ (k\rightarrow \infty)
\end{equation}
for any integer $p\geqslant 1$.
Let $\psi \in C(\overline{\mathbb B}_n)$. It follows from  Lemma \ref{lem:compactness-radial-errors} that
$$
s_k(T_\rho T_{\psi}-T_{\psi}T_\rho)=o(k^{-\frac{1}{n}})\ \ \ \  (k\rightarrow \infty).
$$
Note that
$$
(T_\rho)^pT_{\psi}-T_{\psi}(T_\rho)^p
=\sum_{\ell=0}^{p-1}(T_\rho)^\ell(T_\rho T_{\psi}-T_{\psi}T_\rho)(T_\rho)^{p-1-\ell},
$$
and each summand contains the commutator $T_\rho T_{\psi}-T_{\psi}T_\rho$ and $(p-1)$ copies of $T_\rho$. Thus, we obtain by
\eqref{eq:s_k product estimate-small} that
$$
s_k\!\left[(T_\rho)^\ell(T_\rho T_{\psi}-T_{\psi}T_\rho)(T_\rho)^{p-1-\ell}\right]=o(k^{-\frac{p}{n}})\ \ \ \ (k\rightarrow \infty).
$$
Applying \eqref{eq:s_k-finite-sum}, we conclude that
\begin{equation}\label{eq:s_k commutator}
s_k\!\big[(T_\rho)^pT_{\psi}-T_{\psi}(T_\rho)^p\big]=o(k^{-\frac{p}{n}}),\qquad p\geqslant1.
\end{equation}

We now verify the initial case $m=1$. If $J=\varnothing$, then $v_1=\varphi_1$ and $r=0$. Letting $E=0$, the result follows immediately. If $J=\{1\}$, then $v_1=\rho \varphi_1$ and $r=1$. Taking
$$E=T_{\rho \varphi_1}-T_{\varphi_1}T_\rho,$$
Lemma \ref{lem:compactness-radial-errors} gives
$$
s_k(T_{\rho \varphi_1}-T_{\varphi_1}T_\rho) = o(k^{-\frac{1}{n}}) \ \ \ \ (k\rightarrow \infty).
$$
This proves the assertion for one factor.

Suppose that $m\geqslant 2$ and the assertion holds for every product of $(m-1)$ factors. Let
$$ p=\bigl|J\cap\{1, 2, \cdots,m-1\}\bigr|. $$
By the induction hypothesis,
\begin{equation}\label{eq:guina-m-1}
T_{v_1}T_{v_2}\cdots T_{v_{m-1}}=T_{\varphi_1\varphi_2\cdots \varphi_{m-1}}(T_\rho)^p+E_0,\qquad s_k(E_0)=o(k^{-\frac{p}{n}}) \ \ \ \ (k\rightarrow \infty),
\end{equation}
where $E_0$ is compact. In order to prove the assertion for \(m\) factors, we consider the following two cases:

\textbf{Case 1.} If $m\notin J$, then $v_{m}=\varphi_{m}$ and $r=|J|=p$. The case $r=0$ is immediate, so we may assume that $r\geqslant1$.  By \eqref{eq:guina-m-1}, we have
$$
\begin{aligned}
&T_{v_1}T_{v_2}\cdots T_{v_{m}}-T_{\varphi_1 \varphi_2\cdots \varphi_{m}}(T_\rho)^r \\
&=E_0T_{\varphi_{m}}+T_{\varphi_1\varphi_2\cdots \varphi_{m-1}}(T_\rho)^rT_{\varphi_{m}}-T_{\varphi_1 \varphi_2\cdots \varphi_{m}}(T_\rho)^r \\
&=E_0T_{\varphi_{m}}+T_{\varphi_1\varphi_2\cdots \varphi_{m-1}}\bigl[(T_\rho)^rT_{\varphi_{m}}-T_{\varphi_{m}}(T_\rho)^r\bigr]+\bigl[T_{\varphi_1\varphi_2\cdots \varphi_{m-1}}T_{\varphi_{m}}-T_{\varphi_1\varphi_2\cdots \varphi_{m}}\bigr](T_\rho)^r\\
&:=\mathrm I+\mathrm{II}+\mathrm{III}.
\end{aligned}
$$
Denote $$E=T_{v_1}T_{v_2}\cdots T_{v_{m}}-T_{\varphi_1 \varphi_2\cdots \varphi_{m}}(T_\rho)^r .$$ Since $T_{\rho}$ and $E_0$ are both compact, we deduce that
$E$ is also compact.  In addition, combining
Lemma \ref{lem:singular-value inequalities} and \eqref{eq:guina-m-1} gives
$$s_k(\mathrm I)=s_k(E_0T_{\varphi_{m}})\leqslant\|\varphi_{m}\|_\infty s_k(E_0)=o(k^{-\frac{r}{n}}).$$
Using  the singular-value estimate of the commutator $(T_\rho)^pT_{\varphi_{m}}-T_{\varphi_{m}}(T_\rho)^p$ in \eqref{eq:s_k commutator},
we have that
$$s_k(\mathrm{II})=s_k\big[T_{\varphi_1\varphi_2\cdots \varphi_{m-1}}\bigl((T_\rho)^rT_{\varphi_{m}}-T_{\varphi_{m}}(T_\rho)^r\bigr)\big]= o(k^{-\frac{r}{n}}),$$
since $T_{\varphi_1 \varphi_2\cdots \varphi_{m-1}}$ is bounded. As $T_{\varphi_1\varphi_2\cdots \varphi_{m-1}}T_{\varphi_{m}}-T_{\varphi_1\varphi_2\cdots \varphi_{m}}$
is compact, using the estimate in \eqref{eq:s_k111111} and $(c)$ in Lemma \ref{lem:singular-value inequalities} we obtain
$$
s_k\!\left[\bigl(T_{\varphi_1\varphi_2\cdots \varphi_{m-1}}T_{\varphi_{m}}-T_{\varphi_1\varphi_2\cdots \varphi_{m}}\bigr)(T_\rho)^r\right]=o(k^{-\frac{r}{n}})
$$
for $r\geqslant 1$.
Combining the estimates for $s_k(\mathrm{I})$, $s_k(\mathrm{II})$ and $s_k(\mathrm{III})$ gives
$$s_k\!(E)=o(k^{-\frac{r}{n}})\ \ \ \ (k\rightarrow \infty).$$

\textbf{Case 2.} If $m\in J$, then $v_{m}=\rho \varphi_{m}$ and $r=|J|=p+1$. As  the case $p=0$ is immediate,  we assume $p\geqslant1$.
In this case, we write
$$
\begin{aligned}
&T_{v_1}T_{v_2}\cdots T_{v_{m}}-T_{\varphi_1\varphi_2\cdots \varphi_{m}}(T_\rho)^r \\
&=E_0T_{\rho \varphi_{m}}+T_{\varphi_1\varphi_2 \cdots \varphi_{m-1}}(T_\rho)^pT_{\rho \varphi_{m}}-T_{\varphi_1\varphi_2\cdots \varphi_{m}}(T_\rho)^r \\
&=E_0T_{\rho \varphi_{m}}+T_{\varphi_1 \varphi_2\cdots \varphi_{m-1}}(T_\rho)^p\bigl(T_{\rho \varphi_{m}}-T_{\varphi_{m}}T_\rho\bigr)\\
&\quad+T_{\varphi_1\varphi_2\cdots \varphi_{m-1}}\bigl((T_\rho)^pT_{\varphi_{m}}-T_{\varphi_{m}}(T_\rho)^p\bigr)T_\rho+\bigl(T_{\varphi_1\varphi_2\cdots \varphi_{m-1}}T_{\varphi_{m}}-T_{\varphi_1\varphi_2\cdots \varphi_{m}}\bigr)(T_\rho)^{r}\\
&:=\mathrm I'+\mathrm{II}'+\mathrm{III}'+\mathrm{IV}'.
\end{aligned}
$$
Define  $$E=T_{v_1}T_{v_2}\cdots T_{v_{m}}-T_{\varphi_1\varphi_2\cdots \varphi_{m}}(T_\rho)^r. $$
Then  $E$ is compact. By Proposition \ref{lem:s_k(T_rho)-eqvi and s_k(T_rho)-leq}, we have
$$
s_k(T_{\rho \varphi_{m}})=O(k^{-\frac{1}{n}})\ \ \ \ (k\rightarrow \infty).
$$
Using \eqref{eq:s_k product estimate-small} and \eqref{eq:guina-m-1}, we obtain that
$$
s_{2k-1}(E_0T_{\rho \varphi_{m}})\leqslant s_k(E_0)s_k(T_{\rho \varphi_{m}})=o(k^{-\frac{r}{n}}),
$$
which implies that
$$
s_k(\mathrm{I'})=o(k^{-\frac{r}{n}}) \ \ \ \  (k\rightarrow \infty).
$$
By Lemma \ref{lem:compactness-radial-errors},
$$
s_k(T_{\rho \varphi_{m}}-T_{\varphi_{m}}T_\rho)=o(k^{-\frac{1}{n}})\ \ \ \  (k\rightarrow \infty).
$$
Combining this  with \eqref{eq:s_k111111}  and \eqref{eq:s_k product estimate-small}, we obtain
$$
s_k\!\left(\mathrm{II'}
\right)
=o(k^{-\frac{r}{n}})\ \ \ \  (k\rightarrow \infty).
$$
Using the same argument as the one used in estimating
 $s_k(\mathrm{I'})$ and $s_k(\mathrm{II'})$, we have
$$
s_k(\mathrm{III'})=o(k^{-\frac{r}{n}})\ \ \ \  (k\rightarrow \infty).
$$
Furthermore, since $T_{\varphi_1\varphi_2\cdots \varphi_{m-1}}T_{\varphi_{m}}-T_{\varphi_1\varphi_2\cdots \varphi_{m}}$ is compact, we conclude by (c) of Lemma \ref{lem:singular-value inequalities} and (\ref{eq:s_k111111}) that
$$
s_k(\mathrm{IV'})=o(k^{-\frac{r}{n}})\ \ \ \  (k\rightarrow \infty).
$$
Combining the estimates for $s_k(\mathrm{I'})$, $s_k(\mathrm{II'})$, $s_k(\mathrm{III'})$
and $s_k(\mathrm{IV'})$
gives
$$s_k\!(E)=o(k^{-\frac{r}{n}}) \ \ \ \ (k\rightarrow \infty).$$

From the discussions in Cases 1 and 2, we  establish the assertion for $m$ factors and complete the induction. This finishes the proof of Proposition \ref{thm:product-reduction}.
\end{proof}

 Let $\psi$ be continuous on $\overline{\mathbb B}_n$ and $r$ be a positive integer.
  The next proposition shows that the compact operator $T_{\psi}(T_{\rho})^r$ cannot have singular values of the form $o(k^{-\frac{r}{n}})$ if $\psi$ is nonzero somewhere on the sphere.

\begin{proposition}\label{lem:lower-bound}
Suppose that $\psi$ is continuous on $\overline{\mathbb B}_n$ and $\psi|_{\mathbb S_n}\not\equiv0$.
For every integer $r\geqslant1$, the operator $T_{\psi}(T_\rho)^r$
is compact and
\begin{equation*}\label{eq:main-term-not-small}
\limsup_{k\rightarrow \infty} k^{\frac{r}{n}} s_k[T_{\psi}(T_\rho)^r]>0.
\end{equation*}
\end{proposition}

\begin{proof}
Denote $K=T_{\psi}(T_\rho)^r$. Since $T_\rho$ is compact
by Lemma~\ref{lem:radial} and $T_{\psi}$ is bounded,
$K$ is compact.

We first compute the average of
$\langle T_{|\psi|^2}e_\alpha,e_\alpha\rangle$
over each homogeneous degree.
Let $h_\ell=\dim (\mathcal H_\ell)$.
Using \eqref{eq:normalized monomial}, we have that
\begin{equation}\label{eq:homogeneous-kernel-diagonal}
\begin{aligned}
\sum_{|\alpha|=\ell}|e_\alpha(z)|^2
&=\frac{(\ell+n)!}{n!}
  \sum_{|\alpha|=\ell}
  \frac{|z^\alpha|^2}{\alpha!}
=\frac{(\ell+n)!}{n!\,\ell!}
  \bigl(|z_1|^2+|z_2|^2+\cdots+|z_n|^2\bigr)^\ell\\
&=\frac{\ell+n}{n}h_\ell|z|^{2\ell}.
\end{aligned}
\end{equation}
Since $e_\alpha\in L_a^2(\mathbb B_n)$, we have
$
\langle T_{|\psi|^2}e_\alpha,e_\alpha\rangle
=
\int_{\mathbb B_n}
|\psi(z)|^2|e_\alpha(z)|^2\,dV(z).
$
Combining
\eqref{eq:homogeneous-kernel-diagonal}
and \eqref{polar} gives that
\begin{equation}\label{eq:boundary-trace-average}
\begin{aligned}
\frac1{h_\ell}\sum_{|\alpha|=\ell}\langle T_{|\psi|^2}e_\alpha,e_\alpha\rangle
&=\frac{\ell+n}{n}\int_{\mathbb B_n}|\psi(z)|^2|z|^{2\ell}\,dV(z)\\
&=2(\ell+n)\int_0^1 t^{2(\ell+n)-1}\left(\int_{\mathbb S_n}|\psi(t\zeta)|^2\,d\sigma(\zeta)\right)dt.
\end{aligned}
\end{equation}

We next evaluate the limit of
\eqref{eq:boundary-trace-average}.
Let
$
c_\psi=
\int_{\mathbb S_n}
|\psi(\zeta)|^2\,d\sigma(\zeta).
$
The continuity of $\psi$ and the assumption
$\psi|_{\mathbb S_n}\not\equiv0$
imply that $c_\psi>0$. Given $\varepsilon>0$, the uniform continuity of $|\psi|^2$
on $\overline{\mathbb B}_n$ allows us to choose
$0<\eta<1$ such that
$$
\bigl||\psi(t\zeta)|^2-|\psi(\zeta)|^2\bigr|
<\varepsilon,
\qquad
1-\eta<t\leqslant1,\quad
\zeta\in \mathbb S_n.
$$
Using
$$
2(\ell+n)\int_0^1t^{2(\ell+n)-1}\,dt=1,
$$
we have
\begin{align*}
\bigg|\frac1{h_\ell}\sum_{|\alpha|=\ell}\langle T_{|\psi|^2}e_\alpha,e_\alpha\rangle-c_\psi\bigg|&
=\bigg|2(\ell+n)\int_0^1t^{2(\ell+n)-1}\int_{\mathbb S_n}\bigl[|\psi(t\zeta)|^2-|\psi(\zeta)|^2\bigr]\,d\sigma(\zeta)\,dt\bigg|\\
&\leqslant2(\ell+n)\int_0^1t^{2(\ell+n)-1}\int_{\mathbb S_n}\bigl||\psi(t\zeta)|^2-|\psi(\zeta)|^2\bigr|\,d\sigma(\zeta)\,dt\\
&=2(\ell+n)\int_0^{1-\eta}t^{2(\ell+n)-1}\int_{\mathbb S_n}\bigl||\psi(t\zeta)|^2-|\psi(\zeta)|^2\bigr|\,d\sigma(\zeta)\,dt\\
&\quad+2(\ell+n)\int_{1-\eta}^1t^{2(\ell+n)-1}\int_{\mathbb S_n}\bigl||\psi(t\zeta)|^2-|\psi(\zeta)|^2\bigr|\,d\sigma(\zeta)\,dt.
\end{align*}
Since $\psi$ is bounded on $\mathbb B_n$, we get that
\begin{align*}\label{eq:trace-average-error}
\begin{split}
\bigg|\frac1{h_\ell}\sum_{|\alpha|=\ell}\langle T_{|\psi|^2}e_\alpha,e_\alpha\rangle-c_\psi\bigg|
&\leqslant4(\ell+n)\|\psi\|_\infty^2\int_0^{1-\eta}t^{2(\ell+n)-1}\,dt+2(\ell+n)\varepsilon\int_{1-\eta}^1t^{2(\ell+n)-1}\,dt\\
&=2\|\psi\|_\infty^2(1-\eta)^{2(\ell+n)}+\varepsilon\bigl[1-(1-\eta)^{2(\ell+n)}\bigr]\\
&<2\|\psi\|_\infty^2(1-\eta)^{2(\ell+n)}+\varepsilon.
\end{split}
\end{align*}
Since $\varepsilon$ is arbitrary and $\eta\in (0, 1)$, we obtain
\begin{equation}\label{eq:trace-average-limit}
\lim_{\ell\to\infty}
\frac1{h_\ell}\sum_{|\alpha|=\ell}
\langle T_{|\psi|^2}e_\alpha,e_\alpha\rangle
=c_\psi.
\end{equation}

Recall that $D_N=\sum\limits_{\ell=N}^{2N}h_\ell$ in \eqref{eq:band-dimension}.
Grouping the sum according to $|\alpha|=\ell$, we obtain
\begin{align*}\label{eq:band-average-error}
\begin{split}
\bigg|\frac1{D_N}\sum_{N\leqslant|\alpha|\leqslant2N}\langle T_{|\psi|^2}e_\alpha,e_\alpha\rangle-c_\psi\bigg|
&=\bigg|\frac1{D_N}\sum_{\ell=N}^{2N}\sum_{|\alpha|=\ell}\langle T_{|\psi|^2}e_\alpha,e_\alpha\rangle-\frac1{D_N}\sum_{\ell=N}^{2N}h_\ell c_\psi\bigg|\\
&=\bigg|\frac1{D_N}\sum_{\ell=N}^{2N}h_\ell\bigg(\frac1{h_\ell}\sum_{|\alpha|=\ell}\langle T_{|\psi|^2}e_\alpha,e_\alpha\rangle-c_\psi\bigg)\bigg|\\
&\leqslant\frac1{D_N}\sum_{\ell=N}^{2N}h_\ell\bigg|\frac1{h_\ell}\sum_{|\alpha|=\ell}\langle T_{|\psi|^2}e_\alpha,e_\alpha\rangle-c_\psi\bigg|\\
&\leqslant\frac{1}{D_N}\sum_{\ell=N}^{2N}h_\ell\sup_{\ell\geqslant N}\bigg|\frac1{h_\ell}\sum_{|\alpha|=\ell}\langle T_{|\psi|^2}e_\alpha,e_\alpha\rangle-c_\psi\bigg|\\
&=\sup_{\ell\geqslant N}\bigg|\frac1{h_\ell}\sum_{|\alpha|=\ell}\langle T_{|\psi|^2}e_\alpha,e_\alpha\rangle-c_\psi\bigg|.
\end{split}
\end{align*}
By \eqref{eq:trace-average-limit}, we conclude that
\begin{equation}\label{eq:band-trace-average-limit}
\lim_{N\to\infty}
\frac1{D_N}
\sum_{N\leqslant|\alpha|\leqslant2N}
\langle T_{|\psi|^2}e_\alpha,e_\alpha\rangle
=c_\psi.
\end{equation}

We now pass from these diagonal entries to $\|T_{\psi}e_\alpha\|_2^2$.  Note that
\begin{equation}\label{eq:band-semicommutator}
\|T_{\psi}e_\alpha\|_2^2=\langle T_{|\psi|^2}e_\alpha,e_\alpha\rangle+\langle (T_{\overline \psi}T_{\psi}-T_{|\psi|^2})e_\alpha,e_\alpha\rangle
\end{equation}
and  $\Pi_Ne_\alpha=e_\alpha$ for  $N\leqslant|\alpha|\leqslant2N$. Then we have
\begin{align*}
\begin{split}
\bigg|\frac1{D_N}\sum_{N\leqslant|\alpha|\leqslant2N}\Big\langle(T_{\overline \psi}T_{\psi}-T_{|\psi|^2})e_\alpha,e_\alpha\Big\rangle\bigg|
&=\bigg|\frac1{D_N}\sum_{N\leqslant|\alpha|\leqslant 2N}\Big\langle(T_{\overline \psi}T_{\psi}-T_{|\psi|^2})\Pi_Ne_\alpha,e_\alpha\Big\rangle\bigg|\\
&\leqslant\frac1{D_N}\sum_{N\leqslant|\alpha|\leqslant2N}\bigg|\Big\langle(T_{\overline \psi}T_{\psi}-T_{|\psi|^2})\Pi_Ne_\alpha,e_\alpha\Big\rangle\bigg|\\
&\leqslant\frac1{D_N}\sum_{N\leqslant|\alpha|\leqslant2N}\left\|(T_{\overline \psi}T_{\psi}-T_{|\psi|^2})\Pi_Ne_\alpha\right\|_2\|e_\alpha\|_2\\
&\leqslant\frac1{D_N}\sum_{N\leqslant|\alpha|\leqslant2N}\left\|(T_{\overline \psi}T_{\psi}-T_{|\psi|^2})\Pi_N\right\|\|e_\alpha\|_2^2\\
&=\|(T_{\overline\psi}T_{\psi}-T_{|\psi|^2})\Pi_N\|.
\end{split}
\end{align*}
Since $T_{\overline \psi}T_{\psi}-T_{|\psi|^2}$ is compact, Lemma \ref{lem:compact-band} gives
$$
\lim_{N\to\infty}
\big\|
(T_{\overline \psi}T_{\psi}-T_{|\psi|^2})\Pi_N
\big\|=0.
$$
Therefore,
\begin{equation}\label{eq:band-compact-trace}
 \lim_{N\to\infty}\bigg|\frac1{D_N}\sum_{N\leqslant|\alpha|\leqslant2N}\left\langle(T_{\overline \psi}T_{\psi}-T_{|\psi|^2})e_\alpha,e_\alpha\right\rangle
\bigg|=0.
\end{equation}
Summing \eqref{eq:band-semicommutator}
and applying \eqref{eq:band-trace-average-limit}
and \eqref{eq:band-compact-trace}, we obtain
\begin{equation}\label{eq:positive-band-mass}
\lim_{N\to\infty}\frac1{D_N}\sum_{N\leqslant|\alpha|\leqslant2N}\|T_{\psi}e_\alpha\|_2^2= c_\psi.
\end{equation}

Next, we are going to compute the sum of the squared singular values of $K\Pi_N$. Since $\Pi_N^*=\Pi_N$, we have
$$
(K\Pi_N)^*(K\Pi_N)
=\Pi_NK^*K\Pi_N
=|K\Pi_N|^2.
$$
Moreover,
$$
\operatorname{ran}(\Pi_NK^*K\Pi_N)
\subset\operatorname{ran}(\Pi_N),
\qquad
\Pi_NK^*K\Pi_N=0
\quad\text{on} \ \ [\operatorname{ran}(\Pi_N)]^\perp.
$$
Thus the eigenvalues of $(\Pi_NK^*K\Pi_N)\big|_{\operatorname{ran}(\Pi_N)}$, counted with multiplicity
and including zeros, are given by
$$
[s_1(K\Pi_N)]^2, [s_2(K\Pi_N)]^2, \cdots,[s_{D_N}(K\Pi_N)]^2.
$$

Observe that
$$
\{e_\alpha: N\leqslant|\alpha|\leqslant2N\}
$$
is an orthonormal basis of $\operatorname{ran}(\Pi_N)$. Moreover,
for each of these basis vectors,
$$
\Pi_Ne_\alpha=e_\alpha,
\qquad
(T_\rho)^re_\alpha
=\frac{1}{(|\alpha|+n+1)^r}e_\alpha.
$$
By Lemma \ref{lem:singular-value-square-sum}, computing the sum of the squared singular values of $K\Pi_N$ gives
\begin{align}\label{eq:band-square-sum}
\begin{split}
\sum_{j=1}^{D_N}[s_j(K\Pi_N)]^2
&=\sum_{N\leqslant|\alpha|\leqslant2N}\left\langle\Pi_NK^*K\Pi_Ne_\alpha,e_\alpha\right\rangle\\
&=\sum_{N\leqslant|\alpha|\leqslant2N}\left\langle K\Pi_Ne_\alpha,K\Pi_Ne_\alpha\right\rangle\\
&=\sum_{N\leqslant|\alpha|\leqslant2N}\|T_{\psi}(T_\rho)^re_\alpha\|_2^2\\
&=\sum_{N\leqslant|\alpha|\leqslant2N}\left\|T_{\psi}\!\left(\frac{e_\alpha}{(|\alpha|+n+1)^r}\right)\right\|_2^2\\
&=\sum_{N\leqslant|\alpha|\leqslant2N}\frac{\|T_{\psi}e_\alpha\|_2^2}{(|\alpha|+n+1)^{2r}}.
\end{split}
\end{align}
By \eqref{eq:positive-band-mass},
for sufficiently large $N$,
$$
\sum_{N\leqslant|\alpha|\leqslant2N}
\|T_{\psi}e_\alpha\|_2^2
\geqslant\frac{c_\psi}{2}D_N.
$$
For $N\leqslant|\alpha|\leqslant2N$ and $N\geqslant1$,
$$
|\alpha|+n+1
\leqslant2N+n+1
\leqslant(n+3)N.
$$
Thus \eqref{eq:band-square-sum} yields
\begin{equation}\label{eq:band-square-lower}
\sum_{j=1}^{D_N}[s_j(K\Pi_N)]^2
\geqslant cD_NN^{-2r},
\qquad
c=\frac{c_\psi}{2(n+3)^{2r}}>0.
\end{equation}
On the other hand, \eqref{eq:R-spectrum} gives
\begin{equation}\label{eq:band-norm-upper}
\begin{aligned}
\|K\Pi_N\|
\leqslant
\|T_{\psi}\|\,\|(T_\rho)^r\Pi_N\|
\leqslant
\frac{\|\psi\|_\infty}{(N+n+1)^r}
\leqslant\|\psi\|_\infty N^{-r}.
\end{aligned}
\end{equation}

To obtain a lower bound for individual singular values,
choose
$$
0<\delta<
\min\left\{\frac1 2,\frac{c}{4\|\psi\|_\infty^2}\right\},
\qquad
q_N=\lfloor\delta D_N\rfloor.
$$
Now using  \eqref{eq:band-dimension} gives
\begin{equation}\label{eq:band-index-scale}
\frac{q_N}{D_N}\rightarrow\delta \ \ \ \ \mathrm{and}\ \ \ \
\frac{q_N}{N^n}
\rightarrow\frac{\delta(2^n-1)}{n!}>0 \ \ \ \ (N\rightarrow \infty).
\end{equation}
For sufficiently large $N$,
$1\leqslant q_N<D_N$ and $q_N\leqslant2\delta D_N$.
Using \eqref{eq:band-norm-upper}, we obtain
\begin{equation}\label{eq:band-first-values}
\begin{aligned}
\sum_{j=1}^{q_N}[s_j(K\Pi_N)]^2
\leqslant q_N\|K\Pi_N\|^2
\leqslant
2\delta\|\psi\|_\infty^2D_NN^{-2r}
\leqslant\frac c2D_NN^{-2r}.
\end{aligned}
\end{equation}
Subtracting \eqref{eq:band-first-values}
from \eqref{eq:band-square-lower} gives
\begin{equation}\label{eq:band-tail-lower}
\sum_{j=q_N+1}^{D_N}[s_j(K\Pi_N)]^2
\geqslant\frac c2D_NN^{-2r}.
\end{equation}

For an operator $F$ with $\operatorname{rank}(F)<j$,
we have $\operatorname{rank}(F\Pi_N)<j$ and
$$
\|K\Pi_N-F\Pi_N\|\leqslant\|K-F\|.
$$
Therefore, it follows from
\eqref{eq:approximation-number} that
\begin{equation*}\label{eq:band-compression-values}
\begin{aligned}
s_j(K\Pi_N)
&\leqslant
\inf
\{\|K\Pi_N-F\Pi_N\|: \operatorname{rank}(F)<j\}\\
&\leqslant
\inf\{\|K-F\|: {\operatorname{rank}(F)<j}\}
=s_j(K).
\end{aligned}
\end{equation*}
For $j>q_N$, the monotonicity of singular values  implies
$$
s_j(K\Pi_N)\leqslant s_j(K)\leqslant s_{q_N}(K).
$$
Hence
$$
\sum_{j=q_N+1}^{D_N}[s_j(K\Pi_N)]^2
\leqslant(D_N-q_N)[s_{q_N}(K)]^2
\leqslant D_N[s_{q_N}(K)]^2.
$$
Combining this estimate with
\eqref{eq:band-tail-lower}, we obtain
\begin{equation*}\label{eq:band-singular-lower}
s_{q_N}(K)\geqslant\sqrt{\frac c2}\,N^{-r}
\end{equation*}
for sufficiently large $N$.
Using \eqref{eq:band-index-scale}, we conclude that
$$
\limsup_{N\to\infty}
q_N^{\frac{r}{n}}s_{q_N}(K)
\geqslant
\sqrt{\frac c2}
\left(\frac{\delta(2^n-1)}{n!}\right)^{\frac{r}{n}}
>0.
$$
Noting that $q_N\to\infty$ as $N\rightarrow \infty$, we derive that
the sequence
$\{k^{\frac{r}{n}}s_k(K)\}_{k=1}^\infty$ cannot converge to zero.
This completes the proof of Proposition \ref{lem:lower-bound}.
\end{proof}

\section{Localization and boundary factorization}\label{sec:local}

This section connects the estimates of Section \ref{sec:operators} with the local behavior of
harmonic symbols near the sphere. Proposition \ref{lem:localization} shows that, after
multiplication by a cutoff Toeplitz operator, changing symbols outside the chosen
neighborhood produces a difference whose singular values decay faster than every
negative power. Lemma \ref{lem:harmonic} factors a harmonic symbol vanishing on a
 boundary open set as $u=\rho\Phi$ , where $\Phi$ is real analytic nearby and does not vanish
identically on any smaller nonempty boundary open set. These results allow
us to extend the local factors and apply Proposition \ref{thm:product-reduction}.

\begin{proposition}\label{lem:localization}
Let $U\subset\mathbb C^n$ be a nonempty open subset. Suppose that
$v_j\in L^\infty(\B)$ and $f_j\in L^\infty(\B)$ satisfy $v_j=f_j$ almost everywhere
on $U\cap\B$ for $1\leqslant j\leqslant m$.
If $\chi\in C_c^\infty(U)$, i.e., $\chi$ is a smooth function with compact support in $U$, then
\begin{equation*}\label{eq:localized-error}
s_k\!\left[
\left(T_{v_1}T_{v_2}\cdots T_{v_m}-T_{f_1}T_{f_2}\cdots T_{f_m}\right)T_\chi
\right]=O(k^{-L})
\qquad(k\to\infty)
\end{equation*}
for every $L>0$.
\end{proposition}

\begin{proof}

We divide the proof into three steps.

\textbf{Step 1. The case of one factor.}

If $m=1$,  then we have
$$
(T_{v_1}-T_{f_1})T_\chi
=
T_{v_1-f_1}T_\chi.
$$
Since $\operatorname{supp}(\chi)$ is a compact subset
of $U$, we can choose an open set $V$ such that
$$
\operatorname{supp}(\chi)
\subset V\subset\overline V\subset U.
$$
 Choose $\eta\in C_c^\infty(U)$ satisfying $0\leqslant\eta\leqslant1$, and $\eta=1$ on $V$. The assumption on $v_1-f_1$ implies
$$
(v_1-f_1)\eta=0
\qquad\text{a.e. on }\mathbb B_n.
$$
Therefore, we have that
\begin{equation}\label{eq:loc-one-factor}
\begin{aligned}
T_{v_1-f_1}T_\chi f
&=
P\bigl[(v_1-f_1)P(\chi f)\bigr]
\\
&=
P\bigl[(v_1-f_1)\eta P(\chi f)\bigr]
+
P\bigl[(v_1-f_1)(1-\eta)P(\chi f)\bigr]
\\
&=
P\bigl[(v_1-f_1)(1-\eta)P(\chi f)\bigr]
\end{aligned}
\end{equation}
for every $f\in L_a^2(\mathbb B_n)$. Next, we analyze the operator
\begin{align}\label{A}
A_{\eta, \chi}g=(1-\eta)P(\chi g)
\end{align}
on $L^2(\mathbb B_n)$.

\textbf{Substep 1.1. Separation of the supports.}

Since $\eta=1$ on $V$, the function $1-\eta$ vanishes
on an open neighborhood of $\operatorname{supp}(\chi)$.
Hence
$\operatorname{supp}(1-\eta)\cap\overline{\mathbb B}_n$ and
$\operatorname{supp}(\chi)\cap\overline{\mathbb B}_n$
are disjoint compact sets.  If either set is empty, then $A_{\eta, \chi}=0$, and hence \eqref{eq:loc-one-factor} gives $T_{v_1-f_1}T_\chi=0$. Thus, we may assume both $\operatorname{supp}(1-\eta)\cap\overline{\mathbb B}_n$ and $\operatorname{supp}(\chi)\cap\overline{\mathbb B}_n$ are nonempty. The distance between these two compact sets is positive, and there exists a constant $a>0$
such that $|z-w|\geqslant a$
for
$$
z\in\operatorname{supp}(1-\eta)\cap\overline{\mathbb B}_n \ \ \ \ \text{and}\ \ \ \ w\in\operatorname{supp}(\chi)\cap\overline{\mathbb B}_n.
$$
For $z,w\in\overline{\mathbb B}_n$, we have
\begin{align*}
|1-\langle z,w\rangle|
\geqslant 1-\operatorname{Re}\langle z,w\rangle
=\frac{2-|z|^2-|w|^2+|z-w|^2}{2}
\geqslant\frac{|z-w|^2}{2}.
\end{align*}
Taking $\delta=\frac{a^2}2$, we obtain
\begin{equation}\label{eq:loc-denominator-lower-bound}
|1-\langle z,w\rangle|\geqslant\delta>0
\end{equation}
for all $z,w$ such that $$
z\in\operatorname{supp}(1-\eta)\cap\overline{\mathbb B}_n \ \ \ \ \text{and}\ \ \ \ w\in\operatorname{supp}(\chi)\cap\overline{\mathbb B}_n.
$$

\textbf{Substep 1.2. Choosing the expansion centers.}

For each point of $\operatorname{supp}(1-\eta)\cap\overline{\mathbb B}_n$, take the open ball centered at that point with radius $\frac{\delta}2$.
These balls cover $\operatorname{supp}(1-\eta)\cap\overline{\mathbb B}_n$.
Choose a finite subcover and denote its centers by
$z_1, z_2, \cdots,z_s$. Then
$$
z_\nu\in\operatorname{supp}(1-\eta)\cap\overline{\mathbb B}_n \qquad (1\leqslant\nu\leqslant s)
$$
and
$$
\operatorname{supp}(1-\eta)\cap\overline{\mathbb B}_n\subset\bigcup_{\nu=1}^{s}\Big\{z\in\mathbb C^n:|z-z_\nu|<\frac{\delta}2\Big\}.
$$
To assign each point to a single center, define
\begin{equation}\label{eq:loc-partition-definition}
E_\nu=\Big\{z\in\operatorname{supp}(1-\eta)\cap\mathbb B_n:
|z-z_\nu|<\frac{\delta}2,|z-z_\mu|\geqslant\frac{\delta}2
\text{ \ for \ every} \ \mu<\nu
\Big\}.
\end{equation}
 Observe that each $E_\nu$ is measurable and contains the points in the $\nu$-th ball that are not contained in any preceding ball.

Let $1\leqslant\nu<\mu\leqslant s$.
If $z\in E_\nu\cap E_\mu$, the definitions of
$E_\nu$ and $E_\mu$ give
$$
|z-z_\nu|<\frac{\delta}{2}\ \ \ \ \text{and} \ \ \ \
|z-z_\nu|\geqslant\frac{\delta}{2},
$$
which is impossible. Hence we have
$
E_\nu\cap E_\mu=\varnothing$ for $\nu\neq\mu$.
For each $z$ in $\operatorname{supp}(1-\eta)\cap\mathbb B_n$, choose the smallest index $\nu$ such that $|z-z_\nu|<\frac{\delta}{2}$. Then $z\in E_\nu$. Consequently,
\begin{equation*}\label{eq:loc-partition-properties}
\bigcup_{\nu=1}^{s}E_\nu
=
\operatorname{supp}(1-\eta)\cap\mathbb B_n.
\end{equation*}

\textbf{Substep 1.3. Expanding the reproducing  kernel.}

Let $z\in E_\nu$ and $w\in\operatorname{supp}(\chi)\cap\mathbb B_n$. By the Cauchy--Schwarz inequality,
\begin{align*}
|\langle z-z_\nu,w\rangle|
&=
\bigg|
\sum_{j=1}^{n}
(z_j-(z_\nu)_j)\overline w_j
\bigg|
\leqslant
\bigg(
\sum_{j=1}^{n}|z_j-(z_\nu)_j|^2
\bigg)^{\frac{1}{2}}
\bigg(
\sum_{j=1}^{n}|w_j|^2
\bigg)^{\frac{1}{2}}
=
|z-z_\nu|\,|w|,
\end{align*}
where $(z_{\nu})_j$ stands for the $j$-th coordinate of $z_{\nu}$.
Using \eqref{eq:loc-denominator-lower-bound} and \eqref{eq:loc-partition-definition}, we obtain
\begin{equation}\label{eq:loc-expansion-ratio}
\begin{aligned}
\left|
\frac{\langle z-z_\nu,w\rangle}{1-\langle z_\nu,w\rangle}
\right|
\leqslant\frac{|z-z_\nu|\,|w|}{\delta}
\leqslant\frac{|z-z_\nu|}{\delta}<\frac12.
\end{aligned}
\end{equation}
Using
\begin{align*}
\frac{1}{(1-t)^{n+1}}=\sum_{\ell=0}^{\infty}\binom{n+\ell}{n}t^\ell\qquad (|t|<1),
\end{align*}
we deduce by \eqref{eq:loc-expansion-ratio} that
\begin{equation}\label{eq:loc-kernel-expansion}
\begin{aligned}
\frac1{(1-\langle z,w\rangle)^{n+1}}
&=\frac1{(1-\langle z_\nu,w\rangle)^{n+1}}\frac1{\left(1-\frac{\langle z-z_\nu,w\rangle}{1-\langle z_\nu,w\rangle}\right)^{n+1}}\\
&=\frac1{(1-\langle z_\nu,w\rangle)^{n+1}}\sum_{\ell=0}^{\infty}\binom{n+\ell}{n}\left(\frac{\langle z-z_\nu,w\rangle}{1-\langle z_\nu,w\rangle}\right)^\ell\\
&=\sum_{\ell=0}^{\infty}\binom{n+\ell}{n}\frac{\langle z-z_\nu,w\rangle^\ell} {(1-\langle z_\nu,w\rangle)^{n+1+\ell}}\\
&=\sum_{\ell=0}^{\infty}\binom{n+\ell}{n}\frac{1} {(1-\langle z_\nu,w\rangle)^{n+1+\ell}}\left(\sum_{j=1}^{n}(z_j-(z_\nu)_j)\overline w_j\right)^\ell\\
&=\sum_{\ell=0}^{\infty}\binom{n+\ell}{n}\frac{1} {(1-\langle z_\nu,w\rangle)^{n+1+\ell}}\sum_{|\alpha|=\ell}\frac{\ell!}{\alpha!}(z-z_\nu)^\alpha\overline w^\alpha\\
&=\sum_{\alpha\in \mathbb N_0^n}\binom{n+|\alpha|}{n}\frac{|\alpha|!}{\alpha!}\frac{(z-z_\nu)^\alpha\overline w^\alpha}{(1-\langle z_\nu,w\rangle)^{n+1+|\alpha|}}.
\end{aligned}
\end{equation}

\medskip
\noindent
\textbf{Substep 1.4. Constructing finite-rank approximations.}

Using the integral formula for the Bergman projection, we obtain that
\begin{equation}\label{eq:loc-projection-integral}
(1-\eta(z))P(\chi g)(z)
=
\int_{\mathbb B_n}
\frac{(1-\eta(z))\chi(w)}
     {(1-\langle z,w\rangle)^{n+1}}
g(w)\,dV(w), \ \ \ \ g\in L^2(\mathbb B_n).
\end{equation}
We truncate \eqref{eq:loc-kernel-expansion} at total degree $(N-1)$ and substitute the resulting finite sum into \eqref{eq:loc-projection-integral}. Define
\begin{equation*}
\begin{aligned}
(F_Ng)(z)
={}&
\sum_{\nu=1}^{s}
\mathbf1_{E_\nu}(z)(1-\eta(z))
\sum_{|\alpha|<N}
\binom{n+|\alpha|}{n}
\frac{|\alpha|!}{\alpha!}
(z-z_\nu)^\alpha
\\
&\qquad\times
\int_{\mathbb B_n\cap\operatorname{supp}(\chi)}
\frac{\chi(w)\overline w^\alpha g(w)}
     {(1-\langle z_\nu,w\rangle)^{n+1+|\alpha|}}
\,dV(w), \ \ \ \  g\in L^2(\mathbb B_n),
\end{aligned}
\end{equation*}
where $\mathbf1_{E_\nu}$ is the characteristic function on $E_\nu$:
$$
\mathbf1_{E_\nu}(z)=
\begin{cases}
1,&z\in E_\nu,\\
0,&z\notin E_\nu.
\end{cases}
$$
Noting that each integral coefficient in the definition of $F_N$ is independent of $z$, the range of $F_N$
is contained in
$$
\operatorname{span}
\big\{\mathbf1_{E_\nu}(z)(1-\eta(z))(z-z_\nu)^\alpha:
1\leqslant\nu\leqslant s,\ |\alpha|<N
\big\}.
$$
Since the number of multi-indices satisfying $|\alpha|<N$ is
$$
\sum_{\ell=0}^{N-1}\binom{\ell+n-1}{n-1}=\binom{N+n-1}{n},
$$
we have
\begin{equation}\label{eq:loc-finite-rank-bound}
\begin{aligned}
\operatorname{rank}(F_N)
\leqslant s\binom{N+n-1}{n}
=\frac{s}{n!}N(N+1)\cdots(N+n-1)\leqslant C N^n
\end{aligned}
\end{equation}
for all $N\geqslant 1$,
where $C$ is a constant independent of $N$.

Fix $z\in E_\nu$. Since  $E_1, E_2, \cdots,E_s$
are pairwise disjoint, we further have
\begin{align}\label{eq:loc-finite-rank-definition}
\begin{split}
(F_Ng)(z)
={}&(1-\eta(z))
\int_{\mathbb B_n\cap\operatorname{supp}(\chi)}
\chi(w)g(w)
\\
&\qquad\times
\sum_{\ell=0}^{N-1}
\binom{n+\ell}{n}
\frac{\langle z-z_\nu,w\rangle^\ell}
{(1-\langle z_\nu,w\rangle)^{n+1+\ell}}
\,dV(w).
\end{split}
\end{align}
By \eqref{eq:loc-denominator-lower-bound}, \eqref{eq:loc-expansion-ratio} and the third equality in \eqref{eq:loc-kernel-expansion}, we have
\begin{align*}
&\left|
\frac1{(1-\langle z,w\rangle)^{n+1}}-\sum_{\ell=0}^{N-1}\binom{n+\ell}{n}\frac{\langle z-z_\nu,w\rangle^\ell}{(1-\langle z_\nu,w\rangle)^{n+1+\ell}}\right|\\
&\leqslant\delta^{-n-1}\sum_{\ell=N}^{\infty}\binom{n+\ell}{n}2^{-\ell}
=\delta^{-n-1}\sum_{\ell=N}^{\infty}\binom{n+\ell}{n}\left(\frac23\right)^\ell\left(\frac34\right)^\ell\\
&\leqslant\delta^{-n-1}\left(\frac34\right)^N\sum_{\ell=0}^{\infty}\binom{n+\ell}{n}\left(\frac23\right)^\ell=3^{n+1}\delta^{-n-1}\left(\frac34\right)^N.
\end{align*}
Then it follows from  \eqref{eq:loc-projection-integral} and \eqref{eq:loc-finite-rank-definition} that $(1-\eta(z))P(\chi g)(z)-(F_Ng)(z)$ can be written as
\begin{align*}
&(1-\eta(z))
\int_{\mathbb B_n\cap\operatorname{supp}\chi}
\chi(w)g(w)
\left[
\frac1{(1-\langle z,w\rangle)^{n+1}}
-
\sum_{\ell=0}^{N-1}
\binom{n+\ell}{n}
\frac{\langle z-z_\nu,w\rangle^\ell}
{(1-\langle z_\nu,w\rangle)^{n+1+\ell}}
\right]dV(w).
\end{align*}
Using $0\leqslant1-\eta\leqslant1$ and the
Cauchy--Schwarz inequality, we therefore obtain
$$
\begin{aligned}
&\left|
(1-\eta(z))P(\chi g)(z)-(F_Ng)(z)
\right|
\\
&\leqslant
3^{n+1}\delta^{-n-1}
\left(\frac34\right)^N
\int_{\mathbb B_n\cap\operatorname{supp}\chi}
|\chi(w)|\,|g(w)|\,dV(w)
\\
&\leqslant
3^{n+1}\delta^{-n-1}\|\chi\|_\infty
\left(\frac34\right)^N
\int_{\mathbb B_n}|g(w)|\,dV(w)
\\
&\leqslant
3^{n+1}\delta^{-n-1}\|\chi\|_\infty
\left(\frac34\right)^N\|g\|_2,
\end{aligned}
$$
where the last inequality is due to
  $V(\mathbb B_n)=1$.

For $z\notin\operatorname{supp}(1-\eta)$,
both $(1-\eta(z))P(\chi g)(z)$ and $(F_Ng)(z)$
are zero.
Since $E_1, E_2, \cdots,E_s$ cover
$\operatorname{supp}(1-\eta)\cap\mathbb B_n$,
the resulting bound holds for every
$z\in\mathbb B_n$, with a positive constant $C$
independent of $N$, $z$, and $g$.
Consequently, we have
$$
\begin{aligned}
\|(1-\eta)P(\chi g)-F_Ng\|_2^2
&=
\int_{\mathbb B_n}
\big|
(1-\eta(z))P(\chi g)(z)-(F_Ng)(z)
\big|^2\,dV(z)
\\
&\leqslant
C^2\left(\frac34\right)^{2N}
\|g\|_2^2,
\end{aligned}
$$
to get that
\begin{align}\label{finite-rank}
\|(1-\eta)P(\chi g)-F_Ng\|_2\leqslant C\left(\frac34\right)^N\|g\|_2.
\end{align}
Moreover,
\begin{align*}
\|F_Ng\|_2
&\leqslant\|(1-\eta)P(\chi g)\|_2+\|F_Ng-(1-\eta)P(\chi g)\|_2\\
&\leqslant\Big[\|\chi\|_\infty+C\Big(\frac34\Big)^N\Big]\|g\|_2,
\end{align*}
which shows that $F_N$ is bounded. We conclude by \eqref{eq:loc-finite-rank-bound} and \eqref{finite-rank} that $A_{\eta,\chi}$ (defined in (\ref{A})) is compact.

\textbf{Substep 1.5. Estimating the singular values.}

Choose $C\geqslant1$ so that
\eqref{eq:loc-finite-rank-bound} and
\eqref{finite-rank} hold.
For sufficiently large $k$, take
$$
N=\Big\lfloor\Big(\frac{k}{2C}\Big)^{\frac1n}\Big\rfloor.
$$
Then
$$
\operatorname{rank}(F_N)\leqslant CN^n\leqslant\frac{k}{2}<k
$$
and
\begin{align*}
s_k(A_{\eta, \chi})
\leqslant
\|A_{\eta, \chi}-F_N\|\leqslant
C\Big(\frac34\Big)^{\big\lfloor(\frac{k}{2C})^{\frac1n}\big\rfloor}.
\end{align*}
Since
$$
L\log (k)+\Big\lfloor\Big(\frac{k}{2C}\Big)^{\frac1n}\Big\rfloor\log\Big(\frac34\Big)
\rightarrow-\infty\qquad (k\rightarrow \infty)
$$
for every $L>0$,
we have that
\begin{equation*}\label{eq:K3/4}
  k^L
\Big(\frac34\Big)^{
\big\lfloor(\frac{k}{2C})^{\frac1n}\big\rfloor}
\rightarrow 0\qquad (k\rightarrow \infty),
\end{equation*}
and hence
\begin{equation}\label{eq:loc-separated-projection-bound}
s_k(A_{\eta, \chi})
=
O(k^{-L})
\qquad(k\to\infty)
\end{equation}
for every $L>0$.
Now it follows from \eqref{eq:loc-one-factor} and (\ref{eq:sv-bounded}) of Lemma \ref{lem:singular-value inequalities} that
\begin{align*}
s_k(T_{v_1-f_1}T_\chi)
&\leqslant
\|v_1-f_1\|_\infty
s_k(A_{\eta, \chi})
=O(k^{-L})
\end{align*}
for every $L>0$.
This proves the case of $m=1$.

\textbf{Step 2. Localizing the product of the remaining factors.}

Let $m\geqslant2$ and assume that the conclusion of the proposition
holds for $(m-1)$ factors and  for arbitrary bounded symbols
satisfying the hypotheses.

Let $V$ be an open subset such that
$$
\operatorname{supp}(\chi)
\subset V\subset\overline V\subset U.
$$
Choose $\eta_1\in C_c^\infty(U)$ satisfying
$$
0\leqslant\eta_1\leqslant1,
\qquad \eta_1=1\quad\text{on }V.
$$
Next, choose an open neighborhood $V_1$ of $\operatorname{supp}(\eta_1)$ such that
$$
\operatorname{supp}(\eta_1)
\subset V_1\subset\overline{V}_1\subset U
$$
and choose $\eta_2\in C_c^\infty(U)$ satisfying
$$
0\leqslant\eta_2\leqslant1,
\qquad \eta_2=1\quad\text{on }V_1.
$$

Apply the induction hypothesis on $V$ to the symbols $\eta_1f_2,f_3,\cdots,f_m$ and $f_2,f_3,\cdots,f_m$. Since the two lists of symbols agree on $V\cap\mathbb B_n$ and $\chi\in C_c^\infty(V)$,
the induction hypothesis gives that
\begin{equation}\label{eq:loc-modified-tail-bound}
 s_k\left[\left(T_{f_2}T_{f_3}\cdots T_{f_m}-T_{\eta_1f_2}T_{f_3}\cdots T_{f_m}\right)T_\chi\right]=O(k^{-L})\qquad(k\to\infty)
\end{equation}
for every $L>0$. When $m=2$, the second product in \eqref{eq:loc-modified-tail-bound} is $T_{\eta_1f_2}$.

For $f\in L_a^2(\mathbb B_n)$, we have
\begin{equation}\label{eq:loc-tail-decomposition}
\begin{aligned}
&(1-\eta_2)T_{f_2}T_{f_3}\cdots T_{f_m}T_\chi f\\
&=(1-\eta_2)\left(T_{f_2}T_{f_3}\cdots T_{f_m}-T_{\eta_1f_2}T_{f_3}\cdots T_{f_m}\right)T_\chi f\\
&\quad+(1-\eta_2)T_{\eta_1f_2}T_{f_3}\cdots T_{f_m}T_\chi f\\
&=(1-\eta_2)\left(T_{f_2}T_{f_3}\cdots T_{f_m}-T_{\eta_1f_2}T_{f_3}\cdots T_{f_m}\right)T_\chi f\\&\quad+(1-\eta_2)P\left(
\eta_1f_2T_{f_3}\cdots T_{f_m}T_\chi f\right).
\end{aligned}
\end{equation}

For the first term $(1-\eta_2)\left(T_{f_2}T_{f_3}\cdots T_{f_m}-T_{\eta_1f_2}T_{f_3}\cdots T_{f_m}\right)T_\chi$, using  (b) of Lemma \ref{lem:singular-value inequalities} and  \eqref{eq:loc-modified-tail-bound}, we obtain
\begin{align*}
&s_k\left[(1-\eta_2)\bigl(T_{f_2}T_{f_3}\cdots T_{f_m}-T_{\eta_1f_2}T_{f_3}\cdots T_{f_m}\bigr)T_\chi\right]\\
&\quad\leqslant\|1-\eta_2\|_\infty s_k\left[\bigl(T_{f_2}T_{f_3}\cdots T_{f_m}-T_{\eta_1f_2}T_{f_3}\cdots T_{f_m}\bigr)T_\chi\right]\\
&\quad\leqslant s_k\left[\bigl(T_{f_2}T_{f_3}\cdots T_{f_m}-T_{\eta_1f_2}T_{f_3}\cdots T_{f_m}\bigr)T_\chi\right]\\
&\quad=O(k^{-L})
\end{align*}
for every $L>0$.

To estimate the singular values of the second term $(1-\eta_2)P\left(
\eta_1f_2T_{f_3}\cdots T_{f_m}T_\chi\right)$,
we write
$$(1-\eta_2)P\left(
\eta_1f_2T_{f_3}\cdots T_{f_m}T_\chi f\right)=A_{\eta_2, \eta_1}Bf, \ \ \ \ f\in L_a^2(\mathbb B_n),$$
where $B: L_a^2(\mathbb B_n)\rightarrow L^2(\mathbb B_n)$ is the operator defined by
 $$Bf=f_2T_{f_3}\cdots T_{f_m}T_{\chi}f, \ \ \ \ f\in L_a^2(\mathbb B_n).$$
For $A_{\eta_2, \eta_1}$,  applying
\eqref{eq:loc-separated-projection-bound}
with $\eta=\eta_2$ and $\chi=\eta_1$ (the required support condition is satisfied because
$\eta_2=1$ on a neighborhood of
$\operatorname{supp}(\eta_1)$) yields
\begin{equation}\label{eq:loc-two-cutoff-bound}
s_k(A_{\eta_2, \eta_1})
=
O(k^{-L})
\qquad(k\to\infty)
\end{equation}
for every $L>0$. For every $f\in L_a^2(\mathbb B_n)$, we have
\begin{align*}
\|Bf\|_2
&\leqslant
\|f_2\|_\infty
\|T_{f_3}\cdots T_{f_m}T_\chi f\|_2
\\
&\leqslant
\|f_2\|_\infty
\bigg(\prod_{j=3}^{m}\|T_{f_j}\|\bigg)
\|T_\chi\|\,\|f\|_2
\\
&\leqslant
\bigg(\prod_{j=2}^{m}\|f_j\|_\infty\bigg)
\|\chi\|_\infty\,\|f\|_2.
\end{align*}
Applying (b) of Lemma \ref{lem:singular-value inequalities}  and \eqref{eq:loc-two-cutoff-bound}, we obtain
\begin{align*}
s_k(A_{\eta_2, \eta_1}B)&\leqslant s_k(A_{\eta_2, \eta_1})\|B\|
\leqslant\bigg(\prod_{j=2}^{m}\|f_j\|_\infty\bigg)\|\chi\|_\infty s_k(A_{\eta_2, \eta_1})
=O(k^{-L})
\end{align*}
for every $L>0$.

Therefore, both operators in \eqref{eq:loc-tail-decomposition} satisfy the required estimate. It follows  from (\ref{finitesum}) that
\begin{equation}\label{eq:loc-tail-bound}
s_k\left[(1-\eta_2)T_{f_2}T_{f_3}\cdots T_{f_m}T_\chi\right]=O(k^{-L})\qquad(k\to\infty)
\end{equation}
for every $L>0$.

\textbf{Step 3. Adding the first factor.}

Since $v_1-f_1=0$ a.e. on
$U\cap\mathbb B_n$ and
$\operatorname{supp}(\eta_2)\subset U$,
$$
(v_1-f_1)\eta_2=0
\quad\text{a.e. on } \  \mathbb B_n.
$$
For every $f\in L_a^2(\mathbb B_n)$, we therefore have
\begin{equation}\label{eq:loc-first-difference-factor}
\begin{aligned}
&T_{v_1-f_1}T_{f_2}\cdots T_{f_m}T_\chi f\\
&=P\left[(v_1-f_1)T_{f_2}\cdots T_{f_m}T_\chi f\right]\\
&=P\left[(v_1-f_1)(1-\eta_2)T_{f_2}\cdots T_{f_m}T_\chi f\right].
\end{aligned}
\end{equation}
Combining \eqref{eq:loc-tail-bound} and \eqref{eq:loc-first-difference-factor}, we obtain
\begin{equation}\label{eq:loc-first-difference-bound}
\begin{aligned}
&s_k\left(T_{v_1-f_1}T_{f_2}\cdots T_{f_m}T_\chi\right)\\
&\quad\leqslant\|v_1-f_1\|_\infty s_k\left((1-\eta_2)T_{f_2}\cdots T_{f_m}T_\chi\right)\\
&\quad=O(k^{-L})\qquad(k\to\infty)
\end{aligned}
\end{equation}
for every $L>0$.

We next apply the induction hypothesis to
$v_2, v_3, \cdots,v_m$ and $f_2, f_3, \cdots, f_m$.
These symbols agree almost everywhere on
$U\cap\mathbb B_n$, so
\begin{equation*}\label{eq:loc-remaining-product-difference}
s_k\left[
\left(
T_{v_2} T_{v_3}\cdots T_{v_m}
-
T_{f_2} T_{f_3}\cdots T_{f_m}
\right)T_\chi
\right]
=
O(k^{-L})
\qquad(k\to\infty)
\end{equation*}
for every $L>0$.
Since $T_{v_1}$ is bounded,
\begin{equation}\label{eq:loc-prefixed-product-difference}
s_k\left[
T_{v_1}
\left(
T_{v_2}T_{v_3}\cdots T_{v_m}
-
T_{f_2}T_{f_3}\cdots T_{f_m}
\right)T_\chi
\right]
=
O(k^{-L})
\qquad(k\to\infty)
\end{equation}
for every $L>0$.

Finally, using
\begin{equation*}\label{eq:loc-product-difference}
\begin{aligned}
&\left(T_{v_1}T_{v_2}\cdots T_{v_m}-T_{f_1}T_{f_2}\cdots T_{f_m}\right)T_\chi\\
&=T_{v_1}\left(T_{v_2}\cdots T_{v_m}-T_{f_2}T_{f_3}\cdots T_{f_m}\right)T_\chi+(T_{v_1}-T_{f_1})T_{f_2}\cdots T_{f_m}T_\chi\\
&=T_{v_1}\left(T_{v_2}\cdots T_{v_m}-T_{f_2}T_{f_3}\cdots T_{f_m}\right)T_\chi+T_{v_1-f_1}T_{f_2}\cdots T_{f_m}T_\chi,
\end{aligned}
\end{equation*}
we have by
\eqref{eq:loc-first-difference-bound} and \eqref{eq:loc-prefixed-product-difference}
that
$$
s_k\left[\left(T_{v_1}\cdots T_{v_m}-T_{f_1}\cdots T_{f_m}\right)T_\chi\right]
=
O(k^{-L})
\qquad(k\to\infty)
$$
for every $L>0$.
This completes the induction and proves Proposition \ref{lem:localization}.
\end{proof}

Next, we  show that a harmonic symbol vanishing on an open subset of the sphere
admits a local factorization $u=\rho\Phi$, where $\rho(z)=1-|z|^2$.
In the proof of the main theorem, we extend these locally defined functions
$\Phi$ to the closed unit ball and use Proposition~\ref{lem:localization}
to localize the operator product.
Then we can  apply Proposition \ref{thm:product-reduction}
to extract the factors $T_\rho$.

\begin{lemma}\label{lem:harmonic}
Let $u$ be a nonzero complex-valued harmonic function on $\mathbb B_n$, continuous on $\overline{\mathbb B}_n$. Suppose that $u=0$ on a nonempty
open subset $W$ of $\mathbb S_n$. For every nonempty open subset $W'$ of $\mathbb S_n$ satisfying $\overline{W'}\subset W$, there exist an open neighborhood $U$ of $\overline{W'}$ and a real-analytic function $\Phi: U\to\mathbb C$ such that
\begin{equation*}\label{eq:harmonic-factorization}
u(z)=\rho(z) \Phi(z), \ \ \ \ z\in U\cap\mathbb B_n.
\end{equation*}
Moreover, $\Phi$ does not vanish identically on any nonempty open subset of $W'$.
\end{lemma}

\begin{proof}
Let
$$
U=\left\{t\zeta:\zeta\in W,\ \frac12<t<2\right\}.
$$
Then $U$ is an open neighborhood of $\overline{W'}$. Moreover, $U$ is symmetric about $\mathbb S_n$ and
$$
U\cap\mathbb S_n=W.
$$
The function $u$ is harmonic on $U\cap\mathbb B_n$, continuous on $U\cap\overline{\mathbb B}_n$, and zero on $U\cap\mathbb S_n$.
Applying  the Schwarz reflection principle (see \cite[Theorem~4.15]{ABR}) to each connected component of $U$, we conclude that $u$ extends harmonically to $U$. Let us  continue to denote this extension by $u$. Then we have
by \cite[Theorem~1.28]{ABR} that $u$ is real analytic on $U$.

We next factor out $\rho$ from this extension. For $\zeta\in W$ and $\frac12<t<2$, using $u(\zeta)=0$ we obtain that
$$
\begin{aligned}
u(t\zeta)
=u(t\zeta)-u(\zeta)
=\int_1^t \frac{d}{dr}u(r\zeta)\,dr
=(t-1)\int_0^1 \left. \frac{d}{dr}u(r\zeta) \right|_{r=1+s(t-1)} \,ds.
\end{aligned}
$$
Define $\Phi$ on $U$ by
$$
\Phi(t\zeta) = -\frac{1}{1+t}\int_0^1\left.\frac{d}{dr}u(r\zeta)\right|_{r=1+s(t-1)}\,ds,\qquad\zeta\in W,\quad \frac12<t<2.
$$

We next prove that $\Phi$ is real analytic on $U$.
 For $\frac12<t<2$ and $0\leqslant s\leqslant1$,  the number $1+s(t-1)$ lies between $1$ and $t$.
Hence
$$
\frac12<1+s(t-1)<2,\qquad[1+s(t-1)]\zeta\in U\quad(\zeta\in W).
$$
For $t\ne1$, the defining formula for $\Phi$ gives
$$
\Phi(t\zeta)=\frac{u(t\zeta)}{1-t^2},
$$
which implies that $\Phi$ is real analytic on $U\backslash \mathbb S_n$.
Fix $\zeta_0\in W$ and choose real-analytic local
coordinates $\zeta=\zeta(y)$ on the sphere, with
$\zeta(y_0)=\zeta_0$.
Since $u$ is real analytic, it has a convergent
expansion near $(1,y_0)$ of the form
$$
u(t\zeta(y))
=
\sum_{q=0}^{\infty}
\sum_{\alpha\in\mathbb N_0^{2n-1}}
c_{q,\alpha}(t-1)^q(y-y_0)^\alpha.
$$
The identity $u(\zeta(y))=0$ implies
$c_{0,\alpha}=0$ for every $\alpha$.
Therefore,
$$
u(t\zeta(y))
=
(t-1)
\sum_{q=1}^{\infty}
\sum_{\alpha\in\mathbb N_0^{2n-1}}
c_{q,\alpha}(t-1)^{q-1}(y-y_0)^\alpha.
$$
It follows that
$$
\Phi(t\zeta(y))
=
-\frac1{1+t}
\sum_{q=1}^{\infty}
\sum_{\alpha\in\mathbb N_0^{2n-1}}
c_{q,\alpha}(t-1)^{q-1}(y-y_0)^\alpha, \ \ \ \ t\neq 1.
$$
The right-hand side is real analytic in a smaller
neighborhood of $(1,y_0)$.
Since $\Phi$ is continuous by its integral formula,
this equality also holds at $t=1$.
The map $(t,y)\mapsto t\zeta(y)$ is a real-analytic
change of coordinates near $\zeta_0$, so $\Phi$
is real analytic near $\zeta_0$.
As $\zeta_0\in W$ is arbitrary, we deduce that $\Phi$ is real
analytic on $U$ and $u=\rho \Phi$ on $U\cap\mathbb B_n$.

It remains to prove that $\Phi$ cannot vanish identically on a nonempty  open subset of $W'$.
Suppose, to the contrary, that $\Phi=0$ on a nonempty open subset $O\subset W'$.  For $\zeta\in O$ and real $t$ sufficiently close to $1$,
$$
u(t\zeta)=(1-t^2)\Phi(t\zeta).
$$
Differentiating with respect to $t$, we obtain
$$
\frac{\partial}{\partial t}u(t\zeta)
=
-2t\,\Phi(t\zeta)
+
(1-t^2)\frac{\partial}{\partial t}\Phi(t\zeta).
$$
Consequently,
\begin{equation}\label{eq:zero-cauchy-data}
u(\zeta)=0,
\qquad
\left.
\frac{\partial}{\partial t}u(t\zeta)
\right|_{t=1}
=
-2\Phi(\zeta)=0,
\ \  \ \ \zeta\in O.
\end{equation}

In the following, we use the harmonicity of $u$ to show that all higher radial derivatives also vanish. The polar-coordinate formula for the Laplacian in Lemma \ref{lem:polar-laplacian} gives
\begin{equation}\label{eq:polar-laplacian}
\frac{\partial^2}{\partial t^2}u(t\zeta)
+
\frac{2n-1}{t}
\frac{\partial}{\partial t}u(t\zeta)
+
\frac{1}{t^2}
\Delta_{\mathbb S_n}u(t\zeta)
=
0
\end{equation}
for $t<1$ sufficiently close to $1$ and $\zeta\in O$. For each fixed $t$, the sphere Laplacian acts on the function $\zeta\mapsto u(t\zeta)$.

We prove by induction that
\begin{equation}\label{eq:all-radial-derivatives-zero}
\left.\frac{\partial^q}{\partial t^q}u(t\zeta)\right|_{t=1}=0, \ \ \ \
q\geqslant0, \quad \zeta\in O.
\end{equation}
The cases $q=0$ and $q=1$ follow from \eqref{eq:zero-cauchy-data} immediately. Suppose that
$$
\left.
\frac{\partial^j}{\partial t^j}u(t\zeta)
\right|_{t=1}
=
0
$$
for  $0\leqslant j\leqslant q+1$  and $\zeta\in O$.  Differentiating \eqref{eq:polar-laplacian}
$q$ times with respect to $t$ and applying
the Leibniz rule, we obtain
\begin{equation}\label{eq:radial-derivative-recursion}
\begin{aligned}
\frac{\partial^{q+2}}{\partial t^{q+2}}u(t\zeta)
={}&
-(2n-1)
\sum_{j=0}^{q}\binom{q}{j}
\left(
\frac{d^{q-j}}{dt^{q-j}}t^{-1}
\right)
\frac{\partial^{j+1}}{\partial t^{j+1}}u(t\zeta)
\\
&-
\sum_{j=0}^{q}\binom{q}{j}
\left(
\frac{d^{q-j}}{dt^{q-j}}t^{-2}
\right)
\Delta_{\mathbb S_n}
\left(
\frac{\partial^j}{\partial t^j}
u(t\zeta)
\right).
\end{aligned}
\end{equation}
At $t=1$, the first sum in
\eqref{eq:radial-derivative-recursion} vanishes
by the induction hypothesis.
For each $0\leqslant j\leqslant q$, the function
$$
\zeta\mapsto
\left.
\frac{\partial^j}{\partial t^j}u(t\zeta)
\right|_{t=1}
$$
is identically zero on $O$.
Its sphere derivatives therefore vanish on $O$, so we have
$$
\left.\Delta_{\mathbb S_n}\left(\frac{\partial^j}{\partial t^j}u(t\zeta)\right)\right|_{t=1}=
\Delta_{\mathbb S_n}\left(\left.\frac{\partial^j}{\partial t^j}u(t\zeta)\right|_{t=1}\right)=0,
\qquad \zeta\in O.
$$
Notice that the second sum in
\eqref{eq:radial-derivative-recursion} also vanishes
at $t=1$. It follows that
$$
\left.
\frac{\partial^{q+2}}{\partial t^{q+2}}u(t\zeta)
\right|_{t=1}
=
0,
\qquad \zeta\in O.
$$
This completes the induction and proves
\eqref{eq:all-radial-derivatives-zero}.

Fix $\zeta_0\in O$ and choose real-analytic local
coordinates
$$
y\longmapsto\zeta(y),
\qquad y\in\mathbb R^{2n-1},
\qquad \zeta(y_0)=\zeta_0,
$$
whose image is contained in $O$. Together with $t$, they give local coordinates $(t,y)\mapsto t\zeta(y)$ near $\zeta_0$. By \eqref{eq:all-radial-derivatives-zero}, we have
$$
\left.\frac{\partial^q}{\partial t^q}u(t\zeta(y))\right|_{t=1}=0
$$
for every $q\geqslant0$ and every $y$  near $y_0$. Differentiating these identities with respect to $y$ gives
\begin{equation}\label{eq:mixed-derivatives-zero}
\left.
\frac{\partial^{|\alpha|}}{\partial y^\alpha}
\frac{\partial^q}{\partial t^q}
u(t\zeta(y))
\right|_{(t,y)=(1,y_0)}
=0,
\qquad
q\geqslant0,\quad
\alpha\in\mathbb N_0^{2n-1}.
\end{equation}

Since the extension $u=\rho \Phi$ and the coordinate map $(t,y)\mapsto t\zeta(y)$ are real analytic, the Taylor
expansion converges near $(1,y_0)$. It follows from  \eqref{eq:mixed-derivatives-zero} that
\begin{equation*}\label{eq:vanishing-taylor-expansion}
\begin{aligned}
u(t\zeta(y))
&=
\sum_{q=0}^{\infty}
\sum_{\alpha\in\mathbb N_0^{2n-1}}
\frac{(t-1)^q(y-y_0)^\alpha}{q!\,\alpha!}
\left.
\frac{\partial^{|\alpha|}}{\partial y^\alpha}
\frac{\partial^q}{\partial t^q}
u(t\zeta(y))
\right|_{(t,y)=(1,y_0)}
=0.
\end{aligned}
\end{equation*}
Consequently, for some sufficiently small
$\varepsilon>0$,
$$
u(t\zeta(y))=0$$
whenever
$1-\varepsilon<t<1$  and $|y-y_0|<\varepsilon$.
Since the set
$$
\bigl\{
t\zeta(y):
1-\varepsilon<t<1,\
|y-y_0|<\varepsilon
\bigr\}
$$
is a nonempty open subset of $\mathbb B_n$, we have by \cite[Theorems 1.27 and 1.28]{ABR} that
$$
u(z)=0, \ \ \ \ z\in \mathbb B_n,
$$
contrary to the hypothesis $u\not\equiv0$.
Thus $\Phi$ does not vanish identically on any nonempty
open subset of $W'$. This finishes the proof of Lemma \ref{lem:harmonic}.
\end{proof}

\section{Proof of the main theorem}\label{sec:main-proof}

In this section, we present the proof of  Theorem \ref{thm:main} by contradiction. The proof uses the local factorization of harmonic symbols near the sphere, together with the singular-value estimates established in Sections \ref{sec:operators} and \ref{sec:local}.

\begin{proof}[Proof of Theorem \ref{thm:main}]
Suppose, to the contrary, that
\begin{equation*}\label{eq:main-contradiction-assumption}
T_{u_1}T_{u_2}\cdots T_{u_m}=0,
\end{equation*}
but $u_j\not\equiv0$ on $\mathbb B_n$ for every $1\leqslant j\leqslant m$. For the sake of clarity, we divide the proof into five steps.

\textbf{Step 1. Showing that $u_1u_2\cdots u_m=0$ on $\mathbb S_n$.}

Since each $u_j$ is continuous on $\overline{\mathbb B}_n$, we have by \cite[Theorem 1]{Coburn1973} that the essential spectra of $T_{u_1}T_{u_2}\cdots T_{u_m}$ and $T_{u_1u_2\cdots u_m}$ satisfy
$$\sigma_{\mathrm e}(T_{u_1}T_{u_2}\cdots T_{u_m})=\sigma_{\mathrm e}(T_{u_1u_2\cdots u_m})
=\{u_1(\zeta)u_2(\zeta)\cdots u_m(\zeta): \zeta \in \mathbb{S}_n\}.$$
It follows that
$$\prod_{j=1}^{m}u_j(\zeta)=0,
\qquad
\zeta\in\mathbb S_n.$$

\textbf{Step 2. Choosing a common boundary set and
obtaining local factorizations.}

We have proved that
$$
u_1(\zeta)u_2(\zeta)\cdots u_m(\zeta)=0
\qquad (\zeta\in\mathbb S_n)
$$
in Step 1. Thus, at each boundary point at least one symbol vanishes,
but the vanishing symbol may depend on the point.
We now choose a common nonempty open subset of the sphere
on which each symbol is either identically zero or nowhere zero.

Let us start with a nonempty open subset of $\mathbb S_n$ and consider $u_1$. If $u_1$ is identically zero on this set, keep the set.
Otherwise, choose a point $\zeta_1$ in the set such that $u_1(\zeta_1)\neq0$. By continuity, there is a nonempty open neighborhood of $\zeta_1$, contained in the current set, on which
$$
|u_1(\zeta)-u_1(\zeta_1)|<\frac{|u_1(\zeta_1)|}{2}.
$$
For every $\zeta$ in this neighborhood,
\begin{align*}
|u_1(\zeta)|
\geqslant |u_1(\zeta_1)|-|u_1(\zeta)-u_1(\zeta_1)|
>\frac{|u_1(\zeta_1)|}{2}>0.
\end{align*}
Hence $u_1$ is either identically zero or nowhere zero on the chosen open set.

Next, we  consider $u_2$ on the chosen set and argue analogously. Passing to a smaller nonempty open subset preserves the property already obtained for $u_1$: a function that is identically zero remains identically zero on a subset, and a function that is nowhere zero remains nowhere zero on a subset.
Continuing successively with $u_3, u_4,\cdots,u_m$, we obtain a nonempty open subset $W\subset\mathbb S_n$ such that, for each $1\leqslant j\leqslant m$,
\begin{equation*}\label{eq:boundary-symbol-alternatives}
u_j|_W\equiv0\qquad\text{or}\qquad u_j(\zeta)\neq0\quad\text{for}\ \zeta\in W.
\end{equation*}
Define
$$
J=\{j:u_j|_W\equiv0\},
\qquad
r=|J|.
$$
If $J$ were empty, all the symbols would be nonzero at every point of $W$, which is a contradiction. Therefore, we have $r\geqslant1$.

Choose a nonempty open subset $W_0$ of the sphere with $\overline{W}_0\subset W$.
For each $j\in J$, Lemma~\ref{lem:harmonic} gives a real-analytic function $\phi_j$, defined in a neighborhood of $\overline{W}_0$, such that $u_j=\rho \phi_j$ inside the ball in that neighborhood. Since $J$ is finite, we may take a common open neighborhood $U_0\subset\mathbb C^n$ of $\overline{W}_0$ on which all these functions are defined and
\begin{equation}\label{eq:local-symbol-factorizations}
u_j=\rho \phi_j \quad\text{on }U_0\cap\mathbb B_n, \qquad j\in J.
\end{equation}
By Lemma \ref{lem:harmonic} again, for each $j\in J$ and every nonempty open subset $O\subset W_0$, there exists $\zeta_0\in O$ such that $\phi_j(\zeta_0)\neq0$.
By continuity, we can choose a nonempty open neighborhood of $\zeta_0$, contained in $O$, on which
$$
|\phi_j(\zeta)-\phi_j(\zeta_0)|<\frac{|\phi_j(\zeta_0)|}{2}.
$$
For every $\zeta$ in this neighborhood,
\begin{align*}
|\phi_j(\zeta)|
\geqslant|\phi_j(\zeta_0)|-|\phi_j(\zeta)-\phi_j(\zeta_0)|>\frac{|\phi_j(\zeta_0)|}{2}>0.
\end{align*}
Thus, every nonempty open subset of $W_0$
contains a smaller nonempty open subset on which
$\phi_j$ is nowhere zero.

Starting with $W_0$, apply this selection successively
to the finitely many indices $j\in J$.
Each restriction preserves the nonvanishing
of the functions considered previously.
We therefore obtain a common nonempty open subset
of $W_0$ on which every $\phi_j$ ($j\in J$)
is nowhere zero.  Choose a nonempty open subset $W_1$ of the sphere
whose closure is contained in this common open set.
Since each $u_j$ ($j\notin J$) is already nowhere zero
on $W_0$, we obtain
\begin{equation}\label{eq:nonzero-local-coefficients}
\overline{W}_1\subset W_0,
\qquad
\begin{cases}
\phi_j(\zeta)\neq0,& j\in J,\\
u_j(\zeta)\neq0,& j\notin J,
\end{cases}
\qquad \zeta\in\overline{W}_1.
\end{equation}

\textbf{Step 3. Extending the local functions to the closed unit  ball.}

For $j\in J$, the function $\phi_j$ is defined on $U_0$.
We extend these functions to the closed unit ball while
preserving their values near $\overline{W}_1$.

Choose an open neighborhood $U_1$ of $\overline{W}_1$
whose closure is a compact subset of $U_0$.
Choose a function $\eta\in C_c^\infty(U_0)$ such that
$0\leqslant\eta\leqslant1$, and $\eta=1$
on a neighborhood of $\overline{U}_1$.
For each $j\in J$, multiply the original function
$\phi_j$ by $\eta$ on $U_0$ and define it to be
zero outside $U_0$.
Thus, the resulting function is given by
$$
\begin{cases}
\eta(z)\phi_j(z),&z\in U_0,\\
0,&z\notin U_0.
\end{cases}
$$
Since the original $\phi_j$ is real analytic on $U_0$, the product $\eta\phi_j$ is smooth and has compact support in $U_0$. Extending this product by zero outside $U_0$ yields a smooth function on $\mathbb C^n$, which we continue to denote by $\phi_j$ for simplicity.

Since $\eta=1$ on $U_1$, the extension agrees with
the original $\phi_j$ there.
Thus \eqref{eq:local-symbol-factorizations} gives
\begin{equation}\label{eq:global-coefficient-properties}
u_j=\rho\phi_j
\quad\text{on }U_1\cap\mathbb B_n,
\qquad j\in J.
\end{equation}
The extension also preserves the values of $\phi_j$
on $\overline{W}_1\subset U_1$.
Hence \eqref{eq:nonzero-local-coefficients} gives
$\phi_j(\zeta)\neq0$ for every
$\zeta\in\overline{W}_1$ and $j\in J$.
For $j\notin J$, set $\phi_j=u_j$ on
$\overline{\mathbb B}_n$.
These functions are continuous on the closed ball.
By \eqref{eq:nonzero-local-coefficients},
they are also nonzero at every point of $\overline{W}_1$.
Consequently, all the functions $\phi_1, \phi_2,\cdots,\phi_m$
are continuous on $\overline{\mathbb B}_n$ and are nowhere
zero on $\overline{W}_1$.

Now we are ready to define symbols on $\overline{\mathbb B}_n$ by
\begin{equation*}\label{eq:modified-global-symbols}
v_j=
\begin{cases}
\rho\phi_j,&j\in J,\\
\phi_j,&j\notin J.
\end{cases}
\end{equation*}
For $j\in J$, \eqref{eq:global-coefficient-properties} gives $v_j=u_j$ on $U_1\cap\mathbb B_n$.
For $j\notin J$, the same equality follows from $\phi_j=u_j$.
Therefore, we obtain that
\begin{equation}\label{eq:modified-symbols-local-agreement}
v_j=u_j
\quad\text{on }U_1\cap\mathbb B_n,
\qquad 1\leqslant j\leqslant m.
\end{equation}

\textbf{Step 4. Applying localization and the product estimate.}

Choose $\chi\in C_c^\infty(U_1)$ such that $0\leqslant\chi\leqslant1$, and $\chi=1$ on a neighborhood of $\overline{W}_1$. It follows that
$$
T_{v_1}T_{v_2}\cdots T_{v_m}T_\chi=\left(T_{v_1}T_{v_2}\cdots T_{v_m}-T_{u_1}T_{u_2}\cdots T_{u_m}\right)T_\chi.
$$
On the set $U_1\cap\mathbb B_n$, by \eqref{eq:modified-symbols-local-agreement} we have $v_j=u_j$. Proposition \ref{lem:localization} therefore shows that $T_{v_1}T_{v_2}\cdots T_{v_m}T_\chi$ is compact and
\begin{equation}\label{eq:A-rapid}
s_k(T_{v_1}T_{v_2}\cdots T_{v_m}T_\chi)=O(k^{-L})
\qquad(k\to\infty)
\end{equation}
for every $L>0$.
We next apply Proposition~\ref{thm:product-reduction}
to the $(m+1)$ factors in $T_{v_1}T_{v_2}\cdots T_{v_m}T_\chi$, taking
$$
\varphi_j=\phi_j
\quad(1\leqslant j\leqslant m),
\qquad
\varphi_{m+1}=\chi,
$$
and retaining the index set $J$.
Since $m+1\notin J$,
the number of factors carrying the designated factor $\rho$ remains $r$.

Let $\psi=\chi \phi_1\phi_2\cdots \phi_m$. Then Proposition \ref{thm:product-reduction} gives
\begin{equation}\label{eq:A-main-term}
T_{v_1}T_{v_2}\cdots T_{v_m}T_\chi=T_\psi(T_\rho)^r+E,\qquad s_k(E)=o(k^{-\frac{r}{n}})\qquad(k\to\infty).
\end{equation}
Note that every $\phi_j$ is nonzero on $\overline{W}_1$, and $\chi=1$ on $\overline{W}_1$. Hence, we have
\begin{equation}\label{eq:main-coefficient-nonzero}
\psi(\zeta)
=
\phi_1(\zeta)\phi_2(\zeta)\cdots \phi_m(\zeta)\neq0,\qquad \zeta\in \overline{W}_1.
\end{equation}
Thus $\psi\in C(\overline{\mathbb B}_n)$ and its
restriction to the sphere is not identically zero.

\textbf{Step 5. Deriving the contradiction.}

Choosing $L>\frac rn$ in \eqref{eq:A-rapid}, we have
$$
k^{\frac{r}{n}}s_k(T_{v_1}T_{v_2}\cdots T_{v_m}T_\chi)
=
O(k^{\frac{r}{n}-L})
\rightarrow0
\qquad(k\rightarrow\infty).
$$
Recall that \eqref{eq:A-main-term}  gives
$$
k^{\frac{r}{n}}s_k(E)\rightarrow0
\qquad(k\rightarrow\infty).
$$
Since $T_\psi(T_\rho)^r=T_{v_1}T_{v_2}\cdots T_{v_m}T_\chi-E$ and $s_k(-E)=s_k(E)$,
we have by
\eqref{eq:sv-sum} of Lemma \ref{lem:singular-value inequalities} that
\begin{equation}\label{eq:main-odd-index-estimate}
k^{\frac{r}{n}}s_{2k-1}\bigl(T_\psi(T_\rho)^r\bigr)\leqslant k^{\frac{r}{n}}s_k(T_{v_1}T_{v_2}\cdots T_{v_m}T_\chi)+k^{\frac{r}{n}}s_k(E)\rightarrow0
\qquad(k\rightarrow\infty).
\end{equation}

To obtain the estimate for every index, let
$l\geqslant1$ and choose $k=\lfloor\frac{l+1}{2}\rfloor$.
Then $2k-1\leqslant l\leqslant2k$. Since singular values are nonincreasing, \eqref{eq:main-odd-index-estimate} gives
\begin{equation*}
\begin{aligned}
l^{\frac{r}{n}}s_l\bigl(T_\psi(T_\rho)^r\bigr)
&\leqslant
(2k)^{\frac{r}{n}}s_{2k-1}\bigl(T_\psi(T_\rho)^r\bigr)
\\
&=2^{\frac{r}{n}}k^{\frac{r}{n}}s_{2k-1}\bigl(T_\psi(T_\rho)^r\bigr) \rightarrow 0\qquad(l\to\infty).
\end{aligned}
\end{equation*}
Consequently,
$$
s_l\bigl(T_\psi(T_\rho)^r\bigr)
=
o(l^{-\frac{r}{n}})
\qquad(l\to\infty).
$$
This contradicts Proposition \ref{lem:lower-bound}, since $r\geqslant1$, $\psi$ is continuous on $\overline{\mathbb B}_n$, and $\psi$ is nonzero on $\overline{W}_1$ by \eqref{eq:main-coefficient-nonzero}. The contradiction tells us that
at least one $u_j$ is identically zero on the unit ball $\mathbb B_n$, finishing the proof of Theorem \ref{thm:main}.
\end{proof}
\vspace{2mm}

\subsection*{Acknowledgment}
This work was supported by the National Natural Science Foundation of China
(No. 12371125), the Chongqing Natural Science Foundation (No. CSTB2024NSCQ-MSX0177) and the Inner Mongolia Autonomous Region Science and Technology Program Project (Nos. 2025KYPT0098 and 2026KYPT0006).

\end{document}